\documentclass[a4paper,twoside,12pt]{article}
\usepackage{amsfonts,amsmath,amssymb,amsthm}
\usepackage{geometry}
\usepackage{sectsty}
\usepackage{lineno}

\newcommand{\pkg}[1]{{\fontseries{b}\selectfont #1}}
\let\proglang=\textsf

\providecommand{\keywords}[1]{{\textit{Keywords:}} #1}

\usepackage[normalem]{ulem}
\usepackage{hyperref}
\hypersetup{
    colorlinks,%
    citecolor=blue,%
    filecolor=black,%
    linkcolor=blue,%
    urlcolor=blue
}
\usepackage{natbib}
\usepackage[normalem]{ulem}
\usepackage[inline]{enumitem}
\makeatletter
\newcommand{\inlineitem}[1][]{%
\ifnum\enit@type=\tw@
    {\descriptionlabel{#1}}
  \hspace{\labelsep}%
\else
  \ifnum\enit@type=\z@
       \refstepcounter{\@listctr}\fi
    \quad\@itemlabel\hspace{\labelsep}%
\fi}
\makeatother

\newtheorem{theorem}{Theorem}[section]
\newtheorem{lemma}{Lemma}[section]

\newtheorem{remark}{Remark}[section]
  
\usepackage{xspace}

\usepackage{bbm}
\usepackage{dsfont}
\newcommand\bbone{\ensuremath{\mathbbm{1}}}

\newcommand{\E}{\mathbf{E}}
\newcommand{\Var}{\mathbf{V}}

\DeclareMathOperator*{\spann}{span}

\newcommand{\bX}{\boldsymbol{X}}
\newcommand{\bx}{\boldsymbol{x}}
\newcommand{\bk}{\boldsymbol{k}}

\usepackage{graphicx,subfigure}

\RequirePackage{authblk}

\title{\textbf{Nonparametric estimation in a regression model with additive and multiplicative noise}}
\author{
{Christophe Chesneau\thanks{Christophe Chesneau\\
\hspace*{1.8em}Universit\'e de Caen - LMNO, France,
E-mail: christophe.chesneau@unicaen.fr
}
, Salima El Kolei\thanks{Salima El Kolei\\
\hspace*{1.8em}CREST, ENSAI, France, 
E-mail: salima.el-kolei@ensai.fr}
, Junke Kou\thanks{Junke Kou\\
\hspace*{1.8em}Guilin University of Electronic Technology, China, 
E-mail: kjkou@guet.edu.cn}
, Fabien Navarro\thanks{Fabien Navarro\\
\hspace*{1.8em}CREST, ENSAI, France, 
E-mail: fabien.navarro@ensai.fr}}
}
\date{\today}

\usepackage{fancyhdr}
\def\spacingset#1{\renewcommand{\baselinestretch}%
{#1}\small\normalsize} \spacingset{1}
\spacingset{1.5}
\begin{document}

\maketitle
\begin{abstract}
In this paper, we consider an unknown functional estimation problem in a general nonparametric regression model with the feature of having both multiplicative and additive noise.We propose two new wavelet estimators in this general context. We prove that they achieve fast convergence rates under the mean integrated square error over Besov spaces. The obtained rates  have the particularity of being established under weak conditions on the model. A numerical study in a context comparable to stochastic frontier estimation (with the difference that the boundary is not necessarily a production function) supports the theory.
\end{abstract}

\noindent\keywords{Nonparametric regression, multiplicative regression models, nonparametric frontier,  rates of convergence, wavelets.}

\section{Introduction}
We consider a nonparametric regression model with both multiplicative and additive noise. It is defined by $n$ random variables $Y_1,\ldots,Y_n$,  where
\begin{equation}\label{model}
Y_i=f(\bX_i)U_i+ V_i, \quad  i\in \{1,\ldots,n\}, 
\end{equation}
$f$ is an unknown regression function defined on a subset $\Delta$ of $\mathbb{R}^d$, with $d\ge 1$, $\bX_1, \ldots,\bX_n$ are $n$ identically distributed random vectors with support on $\Delta$,  $U_1, \ldots,U_n$ are $n$ identically distributed random variables and $V_1, \ldots,V_n$ are $n$ identically distributed random variables. Moreover, it is supposed that $\bX_i$ and $U_i$ are independent, and $U_i$ and $V_i$ are independent for any $i \in \lbrace 1, \ldots, n \rbrace$. We are interested in the estimation of the unknown function $r:=f^2$ from $(\bX_1,Y_1),\ldots,(\bX_n, Y_n)$; the random vectors $(U_1,V_1),\ldots,(U_n, V_n)$ form the multiplicative-additive noise.
We consider the general formulation of model given by \eqref{model} since besides the theoretical interest it embodies several potential applications. For example, for $U_i=1$, \eqref{model} becomes  the standard nonparametric regression model with additive noise. It has been studied in many papers via various nonparametric methods, including kernel, splines, projection and wavelets methods. See, for instance, the books of \cite{hardle}, \cite{tsybakov} and \cite{comte}, and the references therein.
For $V_i=0$,  \eqref{model} becomes the standard nonparametric regression model with multiplicative noise. Recent studies can be found in \cite{chi,comte} and the references therein.
For $V_i \neq 0$ with the same variance across $i$ a first study, based on a linear wavelet estimator, was proposed by \cite{chesneau2018linear}. In the case where $V_i$ is a function of $\bX_i$, \eqref{model} becomes the popular nonparametric regression model with multiplicative heteroscedastic noises:
\begin{equation*}
Y_i= g(X_i)+f(X_i) U_i
\end{equation*}

In particular, this model is widely used in financial applications, where the aim is to estimate the variance function $r:=f^2$ from the returns of an asset, for instance, to establish confidence intervals/bands for the mean function $g$. Variance estimation is a fundamental statistical problem with wide applications (see \cite{muller1987estimation,hall1990variance,hardletsy,wang, brown, cai} for fixed design and  recently \cite{kulik2011nonparametric,verzelen2018adaptive,shen2019optimal} for random design).

This multiplicative regression model is also popular in various application areas. For example, in econometrics, within deterministic ($V_i=0$) and stochastic ($V_i\neq 0$) non-parametric frontier models. These models can be interpreted as a special case of the model \eqref{model}, where the random variable $U_i$ represents the technical inefficiency of the company and $V_i$ represents noise that disrupts its performance, the nature of which comes from  unanticipated events such that machine failure, strikes, staff strikes, etc. Under monotonicity and concavity assumptions, the regression function $r$ can be viewed in this case as a function of the production set of a firm and its estimation is therefore of paramount importance in production econometrics. Specific estimation methods have been developed, see for instance \cite{farrell1957measurement,de1984measuring,gijbels1999estimation,daouia2005robust}  for deterministic frontiers models and \cite{fan1996semiparametric,kumbhakar2007nonparametric,simar2011stochastic} for stochastic frontier models. For general regression function and general nonparametric setting, we refer to \cite{girard2008frontier,girard2013frontier} and \cite{jirak2014adaptive} for the definitions and properties of robust estimators.

Applications also exist in signal and image processing (\emph{e.g.}, for Global Positioning System
 signal detection \cite{huang} as well as in speckle noise reduction encounter in particular in synthetic-aperture radar images \cite{kuan1985adaptive} or in medical ultrasound images \cite{rabbani2008speckle,mateo2009finding}), where noise sources can be both additive and multiplicative. In this context, one can also cite \cite{korostelev2012minimax} where the author deals with the estimation of the function's support.

The aim of this paper is to develop wavelet methods for the general model \eqref{model}, with a special focus on mild assumptions on the distributions of $U_1$ and $V_1$ (moments of order $4$ will be required, including Gaussian noise). Wavelet methods are of interest in nonparametric statistics thanks to their ability to estimate efficiently a wide variety of unknown functions, including those with spikes and bumps. We refer to \cite{anto2} and \cite{hardle}, and the references therein.  To the best of our knowledge, their development for \eqref{model} taking in full generality is new in the literature. First of all, we construct a linear wavelet estimator using projections of wavelet coefficients estimators.  We evaluate its rate of convergence under the mean integrated square error (MISE) under mild assumptions on the smoothness of $r$; it is assumed that $r$ belongs to Besov spaces. The linear wavelet estimator has the advantage to be simple, but the knowledge of the smoothness of $r$ is necessary to calibrate a tuning parameter which plays a crucial role in the determination of fast rates of convergence. For this reason, an alternative is given by a nonlinear wavelet estimator.  Using a thresholding rule of  wavelet coefficients estimators, we develop a nonlinear wavelet estimator. To reach the goal of mildness assumptions on the model, we use a truncation rule in the definition of the  wavelet coefficients estimators. This technique was introduced by \cite{delyon} in the nonparametric regression estimation setting, and recently improved in \cite{chesneau} (in a multidimensional regression function under mixing dependence framework) and \cite{chaubey} (for a density estimation under a multiplicative censoring problem). The construction of the hard wavelet estimator does not depend on the smoothness of $r$ and we prove that, from a global point of view, it achieves a better rate of convergence under the MISE. In practice, the empirical performance of the estimators developed in this paper depends on the choice of several parameters, the truncation level of the linear estimator as well as the threshold parameter of the non-linear estimator. We propose here a method of automatic selection of these two parameters based on the 2-fold cross-validation method (2FCV) introduced by \cite{nason:96}. A numerical study, in a context similar to stochastic frontier estimation, is being carried out to demonstrate the applicability of this approach.  
 
The rest of the paper is organized as follows. Preliminaries on wavelets are described in Section~\ref{defff}. Section~\ref{estim} specifies some assumptions on the model, presents our wavelet estimators and the main results on their performances. Numerical experiments are presented in Section \ref{simu}. Section~\ref{proof} is devoted to the proofs of the main result. 
 
\section{Preliminaries on wavelets}\label{defff}
\subsection{Wavelet bases on \texorpdfstring{$[0,1]^d$}{TEXT}}
\label{wav}
We begin with a classical notation in wavelet analysis. A multiresolution analysis (MRA) is a sequence of closed subspaces $\{V_{j}\}_{j\in\mathbb{Z}}$ of the square integrable function space $L^{2}(\mathbb{R})$ satisfying the following properties:
\begin{enumerate}[label=(\roman*)]
\item $V_{j}\subseteq V_{j+1}$, $j\in\mathbb{Z}$. $\mathbb{Z}$ denotes the integer set and $\mathbb{N}:=\{n\in\mathbb{Z}, n\geq0\};$
\item $\overline{\bigcup\limits_{j\in\mathbb{Z}} V_{j}}=L^{2}(\mathbb{R})$ (the space $\bigcup\limits_{j\in\mathbb{Z}} V_{j}$ is dense in $L^{2}(\mathbb{R})$);
\item $f(2\cdot)\in V_{j+1}$ if and only if $f(\cdot)\in V_{j}$ for each $j\in\mathbb{Z}$;
\item There exists $\phi \in L^{2}(\mathbb{R})$ (scaling function) such that $\{\phi(\cdot-k), \ k\in\mathbb{Z} \}$ forms an orthonormal basis of $V_{0}=\overline{\spann}\{\phi(\cdot-k)\}$.
\end{enumerate}
See~\cite{meyer} for further details. For the purpose of this paper, we use the compactly supported scaling function $\phi$ of the Daubechies family, and the associated compactly supported wavelet function $\psi$ (see~\cite{daub}). Then we consider the wavelet tensor product bases on $[0,1]^d$ as described in \cite{cohen}. The main lines and notations are described below.  We set  $\Phi(\bx)=\prod_{v=1}^d \phi (x_v)$ and wavelet functions:
$\Psi_u(\bx) =\psi(x_{u})\prod_{\underset{v\not = u}{v=1}}^d\phi(x_{v})$ when $u\in \{1,\ldots,d\}$, and $\Psi_u(\bx) =\prod_{v\in A_u}\psi(x_{v})\prod_{v\not \in A_u}\phi(x_{v}) $ when $u\in \{d+1,\ldots,2^d-1\}$, 
where $(A_u)_{u\in \{d+1,\ldots,2^d-1\}}$ forms the set of all the non-void subsets of $\{1,\ldots,d\}$ of cardinal superior or equal to $2$. For any integer $j$ and any $\bk=(k_1,\ldots,k_d)$, we set $\Phi_{j,\bk}(\bx)=2^{jd/2}\Phi(2^jx_1-k_1, \ldots,2^jx_d-k_d)$, for any $u\in \{1,\ldots,2^d-1\}$,
$\Psi_{j,\bk,u}(\bx)=2^{jd/2}\Psi_{u}(2^jx_1-k_1, \ldots,2^jx_d-k_d)$. 
Now, let us set $\Lambda_j=\{0,\ldots,2^j-1\}^d$.
Then, with an appropriate treatment on the elements which step on the boundaries $0$ and $1$, there exists an integer $\tau$ such that the system $\mathcal{S}=\{\Phi_{\tau,\bk},
\bk \in \Lambda_{\tau}; \ (\Psi_{j,\bk,u})_{u\in \{1,\ldots,2^d-1\}}, \ j\ge \tau , \
\bk\in  \Lambda_{j}\}$ forms an orthonormal basis of $L^2(\lbrack
0,1 \rbrack^d)$. For any integer $j_* \ge \tau$, a function $h \in L^2(\lbrack
0,1 \rbrack^d)$ can be expressed via $\mathcal{S}$ by the following wavelet series:
\begin{equation}\label{mans}
h(\bx)= \sum_{\bk\in  \Lambda_{j_*}}\alpha_{j_*,\bk}\Phi_{j_*,\bk}(\bx)  +\sum_{j= j_*}^{\infty}\sum_{u=1}^{2^d-1}  \sum_{\bk\in  \Lambda_j}\beta_{j,\bk,u}\Psi_{j,\bk,u}(\bx), \quad \bx\in [0,1]^d,
\end{equation}
where
$\alpha_{j,\bk}=\langle h,\Phi_{j,\bk}\rangle_{[0,1]^d}$ and $\beta_{j,\bk,u}=\langle h,\Psi_{j,\bk,u}\rangle_{[0,1]^d}$. 

Also, let us mention that, by construction, $\int_{[0,1]^d} \Phi_{j,\bk}(\bx)d\bx=2^{-jd/2}$ and $\int_{[0,1]^d} \Psi_{j,\bk,u}(\bx)d\bx=0$. 

Let $P_{j}$ be the orthogonal projection operator from $L^{2}([0,1]^{d})$ onto the space $V_{j}$ with the orthonormal basis $\{\Phi_{j,\bk}(\cdot)=2^{jd/2}\Phi(2^{j}\cdot-\bk),\bk\in\Lambda_{j}\}$. Then, for any $h\in L^{2}([0,1]^{d})$,
\begin{equation*}
P_{j}h(\bx)=\sum_{\bk\in\Lambda_{j}}\alpha_{j,\bk}\Phi_{j,\bk}(\bx), \quad \bx\in [0,1]^d.
\end{equation*}

\subsection{Besov spaces}
Besov spaces are important in theory and applications.  They have the features to a wide variety of function spaces as H\"{o}lder and $L^{2}$ Sobolev spaces. Definitions of those spaces are given below.  Suppose that $\phi$ is $m$ regular (\emph{i.e.}, $\phi\in C^{m}$ and $|D^{\alpha}\phi(x)|\leq c(1+|x|^{2})^{-l}$ for each $l\in\mathbb{Z}$, with $\alpha=0,1,\ldots,m$) and consider the wavelet framework defined in Subsection \ref{wav}. Let $h\in L^{p}([0,1]^{d})$,  $p,q\in[1,\infty]$ and $0<s<m$. Then the following assertions are equivalent:

(1)~$h\in B^{s}_{p,q}([0,1]^{d})$; \quad (2)~$\left\lbrace 2^{js}\|P_{j+1}h-P_{j}h\|_{p}\right\rbrace\in l_{q};$ \quad  (3)~$\{2^{j(s-\frac{d}{p}+\frac{d}{2})}\|\beta_{j,.,.}\|_{p}\}\in l_{q}.$\\
The Besov norm of $h$ can be defined by
\begin{equation}
 \|f\|_{B^{s}_{p,q}}:=\|(\alpha_{\tau,.})\|_{p}+\|(2^{j(s-\frac{d}{p}+\frac{d}{2})}\|\beta_{j,.,.}\|_{p})_{j\geq \tau}\|_{q}, ~\text{where}~ \|\beta_{j,.,.}\|_{p}^{p}=\sum\limits_{u=1}^{2^{d}-1}\sum\limits_{\bk\in\Lambda_j}|\beta_{j,\bk,u}|^{p}.\notag
\end{equation}
Further details on Besov spaces are given in \cite{meyer}, \cite{triebel} and \cite{hardle}. 

\section{Assumptions, estimators and main result}\label{estim}
We consider the model \eqref{model} with $\Delta=[0,1]^d$ for the sake of simplicity. Additional technical assumptions are formulated below. 
\begin{enumerate}[label=\textbf{A.\arabic*}, series = quest]
\item \label{hyp:A1} We suppose that $f: [0,1]^d \rightarrow \mathbb{R}$ is bounded from above. 
\item \label{hyp:A2} We suppose that $\bX_1\sim \mathcal{U}([0,1]^d)$.
\item \label{hyp:A3} We suppose that $U_1$ is reduced (mainly for the sake of simplicity in exposition) and has a moment of order $4$. 
\item \label{hyp:A4} We suppose that $V_1$ has a moment of order $4$. 
\end{enumerate}
The two following assumptions involving $V_i$ and $\bX_i$ are complementary and will be considered separately in the study:
\begin{enumerate}[label=\textbf{A.\arabic*},resume* = quest]
\item \label{hyp:A5} We suppose that $\bX_i$ and $V_i$ are independent for any $i \in \lbrace 1, \ldots, n \rbrace$, and $U_1$ is centered or $V_1$ is centered.
\item \label{hyp:A6} We suppose that $V_i=g(\bX_i)$ where $g: [0,1]^d \rightarrow \mathbb{R}$ is known and bounded from above, and $U_1$ is centered. 
\end{enumerate}
These assumptions will be discussed later; some of them can be relaxed. In our main results, we will consider the two following sets of assumptions: 
\begin{enumerate}[label=\textbf{H.\arabic*}]
\item \label{hyp:H1}$=\{\ref{hyp:A1}, \ref{hyp:A2}, \ref{hyp:A3}, \ref{hyp:A4}, \ref{hyp:A5}\}$, \inlineitem \label{hyp:H2}$=\{\ref{hyp:A1}, \ref{hyp:A2}, \ref{hyp:A3}, \ref{hyp:A4}, \ref{hyp:A6}\}$.
\end{enumerate}
As usual in wavelet methods, the first step towards the estimation of $r$ is to consider its wavelet series given by \eqref{mans}. Then we aim to estimate the unknown wavelet coefficients $\alpha_{j,\bk}=\langle r,\Phi_{j,\bk}\rangle_{[0,1]^d}$ and  $\beta_{j,\bk,u}=\langle r,\Psi_{j,\bk,u}\rangle_{[0,1]^d}$ by efficient estimators. In this study, we propose to estimate $\alpha_{j,\bk}$ by 
\begin{equation}\label{ru}
\hat{\alpha}_{j,\bk}:=\frac{1}{n}\sum_{i=1}^{n}Y_i^2\Phi_{j,\bk}(\bX_i) - v_{j,\bk}, 
\end{equation}
where 
\[
v_{j,\bk}:=\left\{
\begin{aligned}
& \E\left[V_1^2\right] 2^{-jd/2}& & \ {\text{under \ref{hyp:A5}}},\\
&  \int_{[0,1]^d}g^2(\bx) \Phi_{j,\bk}(\bx) d\bx &  & \ {\text{under \ref{hyp:A6}}}. 
\end{aligned}
\right.
\]
This is an unbiased estimator of $\alpha_{j,\bk}$ and it converges to  $\alpha_{j,\bk}$  in $L^2$ (see Lemmas \ref{unbiased} and \ref{var2} in Section \ref{proof}). 
On the other side, we propose to estimate $\beta_{j,\bk,u}$ by 
\begin{equation}\label{ru0}
\hat{\beta}_{j,\bk,u}:=\frac{1}{n}\sum_{i=1}^{n}\left(Y_i^2\Psi_{j,\bk,u}(\bX_i) - w_{j,\bk,u}\right)\bbone_{\left\lbrace \left|Y_i^2\Psi_{j,\bk,u}(\bX_i) - w_{j,\bk,u}\right|\leq \rho_n\right\rbrace },
\end{equation}
where $\bbone_A$ denotes the indicator function over an event $A$, $\rho_n:=\sqrt{n/\ln n}$ and 
\[
w_{j,\bk,u}:=\left\{
\begin{aligned}
& 0& & \ {\text{under \ref{hyp:A5}}},\\
&  \int_{[0,1]^d}g^2(\bx) \Psi_{j,\bk,u}(\bx) d\bx &  & \ {\text{under \ref{hyp:A6}}}. 
\end{aligned}
\right.
\]
Due to the thresholding in its definition, this estimator is not unbiased of $\beta_{j,\bk,u}$ but it converges  to $\beta_{j,\bk,u}$ in $L^2$  (see Lemmas \ref{unbiased} and \ref{var2} in Section \ref{proof}). The role of the thresholding is to relax assumptions on $U_1$ and $V_1$; note that only moments of order $4$ is required in \ref{hyp:A3} and \ref{hyp:A4} including uniform or Gaussian distribution. This selection rule has been introduced in a wavelet setting in \cite{delyon}.  It has been recently improved in \cite{chesneau} (in a multidimensional regression function under mixing dependence framework) and \cite{chaubey} (in a density estimation under multiplicative censoring setting). In this study, we adapt it to the general nonparametric regression model \eqref{model}. 

The next step in the construction of our wavelet estimators for $r$ is to expand the most informative of the wavelet coefficients estimators using the initial wavelet basis. We then define the linear wavelet estimator by
\begin{equation*}
\hat{r}^{\mathrm{lin}}_{n}(\bx):=\sum\limits_{\bk\in  \Lambda_{j_*}}\hat{\alpha}_{j_*,\bk}\Phi_{j_*,\bk}(\bx), \quad \bx\in [0,1]^d.
\end{equation*}
We thus have projected the $\hat{\alpha}_{j,\bk}$'s on the father wavelet basis at a certain level $j_*$. 
Despite the simplicity of its construction,  this estimator has a serious drawback:  its performance highly depends on the choice of the level $j_*$. A suitable choice of $j_*$, but depending on the smoothness of $r$, will be specified in our main result. To address this problem, an alternative is proposed by using a hard-thresholding rule that performs a term-by-term selection of the wavelet coefficient estimators $\hat{\beta}_{j,\bk,u}$ and to project them on the original wavelet basis. We define the nonlinear wavelet estimator by
\begin{equation*}
\hat{r}^{\mathrm{non}}_{n}(\bx):=\sum\limits_{\bk\in  \Lambda_{j_*}}\hat{\alpha}_{j_*,\bk}\Phi_{j_*,\bk}(\bx)+\sum_{j= j_*}^{j_{1}}\sum_{u=1}^{2^d-1}  \sum_{\bk\in  \Lambda_j}\hat{\beta}_{j,\bk,u}\bbone_{\{|\hat{\beta}_{j,\bk,u}|\geq\kappa t_{n}\}}\Psi_{j,\bk,u}(\bx), \  \bx\in [0,1]^d,
\end{equation*}
where $t_{n}:=\sqrt{\ln n /n}=\rho_n^{-1}$. The positive integer $j_{1}$ is specified in our main result, while the constant $\kappa$ will be chosen in its proof (see the proof of Lemma \ref{var3}). The idea of keeping the estimators $\hat{\beta}_{j,\bk,u}$ with magnitude greater to $t_n$ is not new; it is a well-known wavelet techniques with strong mathematical and practical results for numerous nonparametric problems; $t_n$ is so-called ``universal threshold''. We refer to \cite{donoho3}, \cite{delyon} and \cite{hardle}. In this study, we describe how to calibrate such estimator when we deal with the general model \eqref{model}. 

In the sequel, we adopt the following notations: $x_{+}:=\max\{x,0\}$. $A\lesssim B$ denotes $A\leq cB$ for some constant $c>0$; $A\gtrsim B$ means $B\lesssim A$; $A\thicksim B$ stands for both $A\lesssim B$ and $B\lesssim A$.

Theorem~\ref{theo} below determines the rates of convergence attained by $\hat{r}^{\mathrm{lin}}_{n}$ and $\hat{r}^{\mathrm{non}}_{n}$ over the MISE.  
\begin{theorem}\label{theo} Consider the problem defined by \eqref{model} under the assumptions \ref{hyp:H1} or \ref{hyp:H2}, let $r\in B^{s}_{p,q}([0,1]^{d})$ with $p,q\in[1,\infty)$, $s>d/p$. Then 
\begin{itemize}
\item the linear wavelet estimator $\hat{r}^{\mathrm{lin}}_{n}$ with $2^{j_{*}}\thicksim n^{\frac{1}{2s'+d}}$ and $s'=s-d(1/p-1/2)_{+}$ satisfies
\begin{subequations}
\begin{equation}
\E\left[\int_{[0,1]^{d}}\left(\hat{r}^{\mathrm{lin}}_{n}(\bx)-r(\bx)\right)^{2}d\bx\right]\lesssim n^{-\frac{2s'}{2s'+d}};
\label{eq11a}
\end{equation}
\item the nonlinear estimator with $2^{j_{*}}\sim n^{\frac{1}{2m+d}}~(m>s)$, $2^{j_{1}}\sim(n/ \ln n)^{\frac{1}{d}}$ satisfies
\begin{equation}
\E\left[\int_{[0,1]^{d}}\left(\hat{r}^{\mathrm{non}}_{n}(\bx)-r(\bx)\right)^{2}d\bx\right]\lesssim(\ln n){n}^{-\frac{2s}{2s+d}}.
\label{eq11b}
\end{equation}
\end{subequations}
\end{itemize}
\end{theorem}
The obtained rates of convergence are those obtained in the standard density estimation problem or the regression function estimation problem under the MISE over Besov spaces (see \cite{hardle}). Under some strong conditions on the model as $U_i:=1$, the rate of convergence $n^{-\frac{2s}{2s+d}}$ is proved to be optimal in the minimax sense (see \cite{hardle} and \cite{tsybakov}). So our nonlinear wavelet estimator can be optimal in the minimax sense up to a $\ln n$. However, in full generality, without specifying the distributions of $U_i$ and $V_i$, the optimal lower bounds for the MISE are difficult to determine via standard techniques (Fano's lemma, \ldots) and the optimality of our estimators remains an open question. 

\begin{remark}
Some assumptions used in Theorem~\ref{theo} can be relaxed without changing the result. In particular, one can consider the domain $\Delta=[a,b]^d$ with $(a,b)\in\mathbb{R}^2$ and $a<b$ with an adaptation of the wavelet basis. In this case, we can also replace \ref{hyp:A2} by $\bX_{1}$ with density function $h: [a,b]^d\rightarrow \mathbb{R}$ bounded from below, with the following wavelet coefficient estimators: 
\begin{equation*}
\hat{\alpha}_{j,\bk}=\frac{1}{n}\sum_{i=1}^{n} \frac{Y_i^2}{h(\bX_i)}\Phi_{j,\bk}(\bX_i) - v_{j,\bk}, 
\end{equation*}
and
\begin{equation*}
\hat{\beta}_{j,\bk,u}=\frac{1}{n}\sum_{i=1}^{n}\left(\frac{Y_i^2}{h(\bX_i)}\Psi_{j,\bk,u}(\bX_i) - w_{j,\bk,u}\right)\bbone_{\left\lbrace \left|\frac{Y_i^2}{h(\bX_i)}\Psi_{j,\bk,u}(\bX_i) - w_{j,\bk,u}\right|\leq\rho_n\right\rbrace}.
\end{equation*}
Finally, note that, in \ref{hyp:A6} can be improved by assuming $g$ unknown. To the best of our knowledge, only \cite{cai} have developed  wavelet methods in this case for $d=1$, deterministic design $(\bX_i:=i/n)$ and infinite moments for $U_i$. Extension of these methods in the general setting of \eqref{model} needs further developments that we leave for a future work.  
\end{remark}

\section{Numerical Experiments}
\label{simu}
To illustrate the empirical performance of the estimators proposed in this work, we carried out a simulation study. The objective is to highlight some of the theoretical findings using numerical examples. We begin by giving some details about the specificities inherent in wavelet estimators in a non-deterministic design framework. We also try to propose a realistic simulation setting using an adaptive selection method  to select both the truncation parameter of the linear estimator and the threshold parameter of the non-linear estimator. In this context, we compare their empirical performances in the model with both multiplicative and additive noise. Simulations were performed using \proglang{R} and in particular the \pkg{rwavelet} package~\cite{chesnav} (available from \url{https://github.com/fabnavarro/rwavelet}). 

\subsection{Computational aspect of wavelets and parameters selection}
For fixed design, thanks to Mallat's pyramidal algorithm \citep{mallat:08}, the computation of wavelet-based estimators is simple and fast. When considering uniform random design, the implementation requires some changes and several strategies have been developed in the literature (see, \emph{e.g.}, \cite{cai:98,hall:97}). For uniform design regression, \cite{cai:99} has proposed  to use an approach in which the wavelet coefficients are computed by a simple application of Mallat's algorithm using the ordered $Y_{i}$'s as input variables. We have followed this approach because it preserves the simplicity of calculation and the efficiency of the equispaced algorithm. In the context of wavelet regression in random design with heteroscedastic noise, \cite{kulik:09} and \cite{navarro17} also adopted this approach. 

\begin{figure}[t]
\centering
\subfigure[\textit{Parabolas}]{\includegraphics[width=0.3\textwidth]{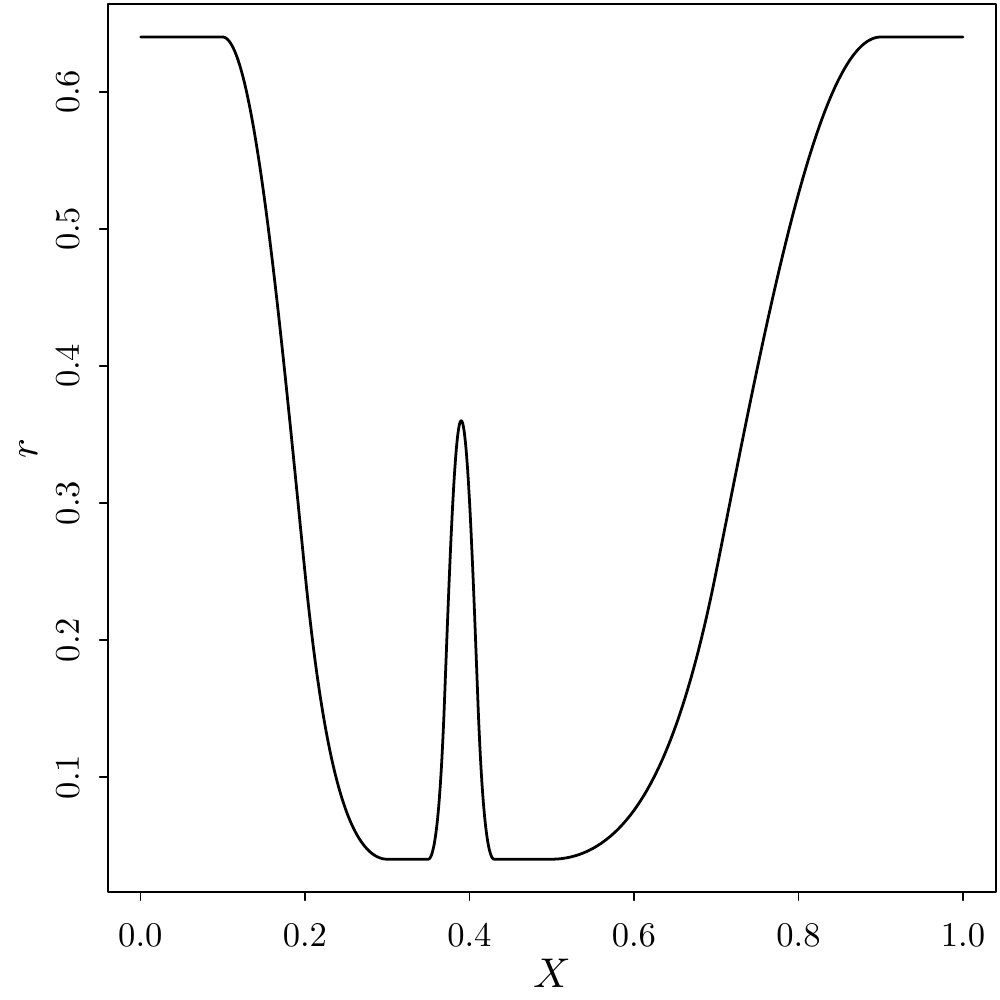}}~
\subfigure[\textit{Ramp}]{\includegraphics[width=0.3\textwidth]{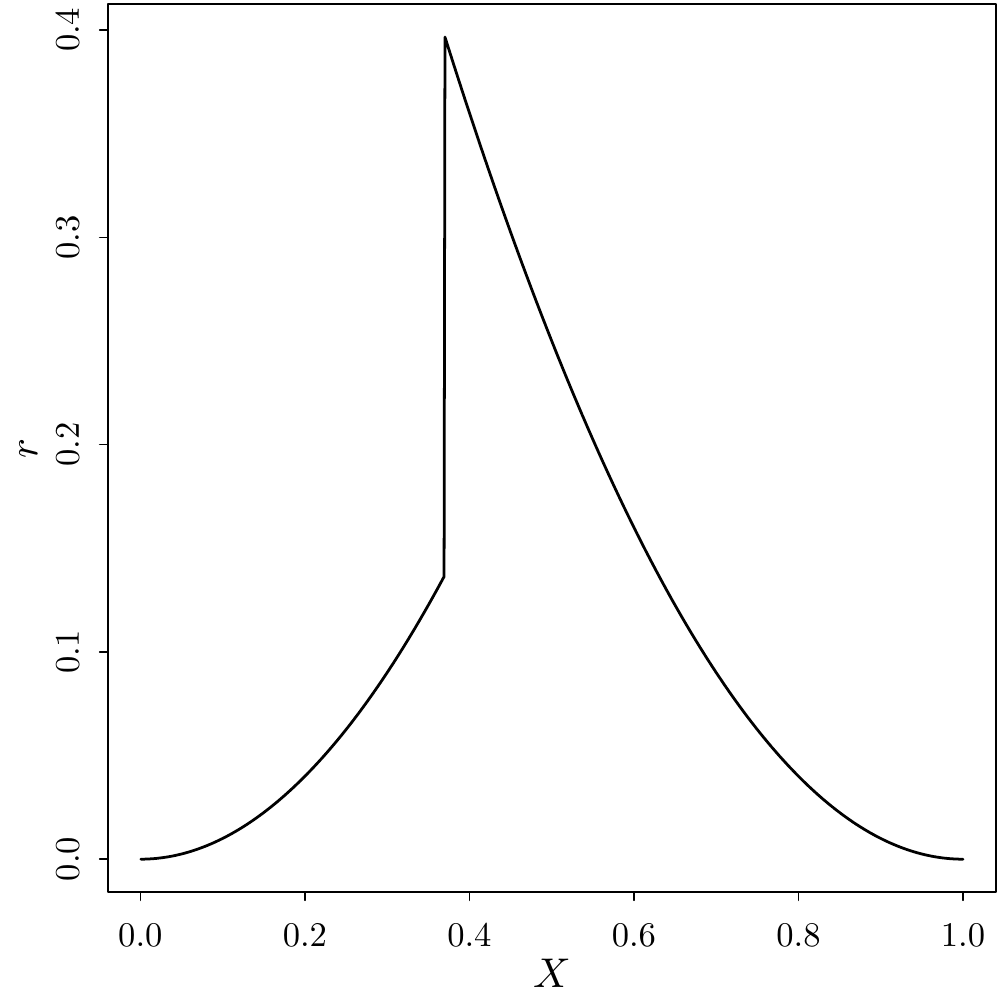}}~
\subfigure[\textit{Blip}]{\includegraphics[width=0.3\textwidth]{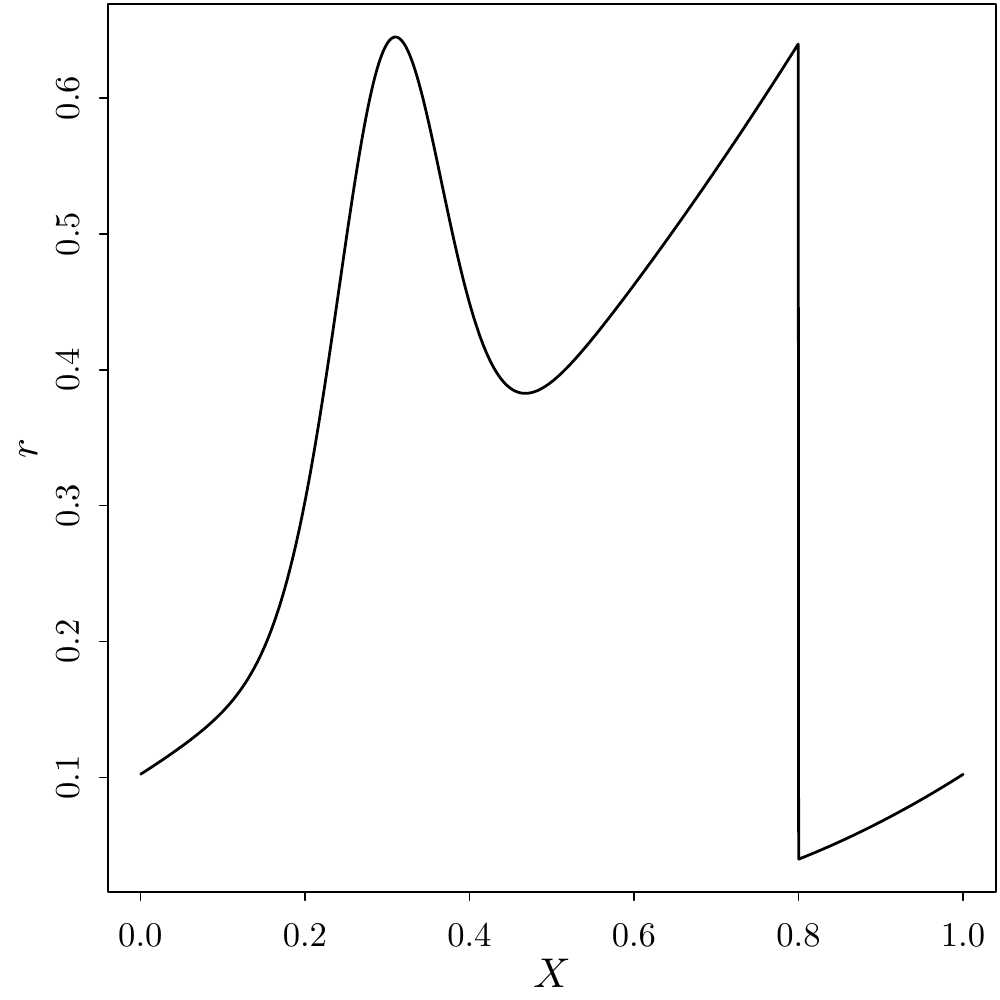}}
\caption{The three test functions considered.}
\label{fig:target}
\end{figure}

Nason successfully adjusted the standard two-Fold Cross Validation (2FCV) method to select the threshold parameter in wavelet shrinkage (see, \cite{nason:96}). For the calibration of linear wavelet estimators, his strategy was used by \cite{navarro17}. We have chosen to apply this approach to select both the threshold and truncation parameter of linear and non-linear estimators. More precisely, in the linear case, we built a collection of linear estimators $\hat{r}_{j_*,n}^{\mathrm{lin}}, j_*=0,1,\ldots,\log2(n)-1$ (by successively adding whole resolution levels of wavelet coefficients), and select the best among this collection by minimizing a 2FCV criterion denoted by $\mathrm{2FCV}(j_*)$. The resulting estimator of the truncation level is denoted by $\hat{j}_*$ and the corresponding estimator of $r$ by $\hat{r}_{\hat{j}_*,n}^{\mathrm{lin}}$ (see, \cite{navarro17,navarro2} for more details). For the nonlinear estimator, the same estimator $\hat{j}_*$ of the truncation parameter $j_*$ obtained for the linear is used. The estimator of the thresholding parameter is obtained using the 2FCV method developed in \cite{nason:96}. The parameter $j_1$ is fixed \textit{a priori} as the maximum level allowed by the wavelet decomposition (\emph{i.e.}, $j_1=\log2(n)-1$). It is a classic choice that allows the coefficients to be selected down to the smallest scale. In addition, in order to facilitate and not to overburden the implementation of the nonlinear estimator, we perform a standard hard thresholding of the wavelet coefficient estimators (rather than the double threshold used in its definition). 
In order to be able to evaluate the performance of these two criteria, the mean square error (MSE) is used (\emph{i.e.}, $\mathrm{MSE}(\hat{r}_{j_*,n},r) = \frac{1}{n}\sum_{i=1}^{n}(r(X_i)-\hat{r}_{j_*,n}(X_i))^2)$). We consider three test functions for $r$ (see  Figure~\ref{fig:target}), commonly used in the wavelet literature, \textit{Parabolas, Ramp} and \textit{Blip} (see, \emph{e.g.}, \cite{donoho3}). 
 In all simulations, we examine the case $d=1$, the design is chosen to satisfy \ref{hyp:A2} (\emph{i.e.}, $\mathcal{U}([0,1])$) and the choice of the wavelet family used is also fixed (\emph{i.e.}, Daubechies compactly supported wavelet with 8 vanishing moments).  

\subsection{Additive-multiplicative regression}
This subsection examines the behaviour and performance of linear and non-linear estimators in the context of additive and multiplicative  regression by considering $V_1\sim \mathcal{N}(0,\sigma^2)$, where $\sigma^2=0.01$ and $U_1\sim U([-1,1])$. Thus, the goal is to estimate the frontier $r$  from $(X_i,Y_i)$ sample simulated from one of the test functions. By applying one of the linear or nonlinear methods developed above to the estimation of $r$, one can construct an estimator whose rate of convergence is given by \eqref{eq11a} and \eqref{eq11b} respectively. Note that here, the nature of the frontier function is not necessarily the same as that commonly found in the literature on stochastic boundary estimation. Indeed, here $r$ is not necessarily a production function (\emph{e.g.}, $r$ is concave), the only assumption we make is given by \ref{hyp:A1}. Thus, the application here can be seen as the estimation of the boundary or frontier of a sample affected by some additive (positive) noise (see \cite{jirak2014adaptive}). 

\begin{figure}[ht]
\centering
\subfigure[]{
\includegraphics[width=0.31\textwidth]{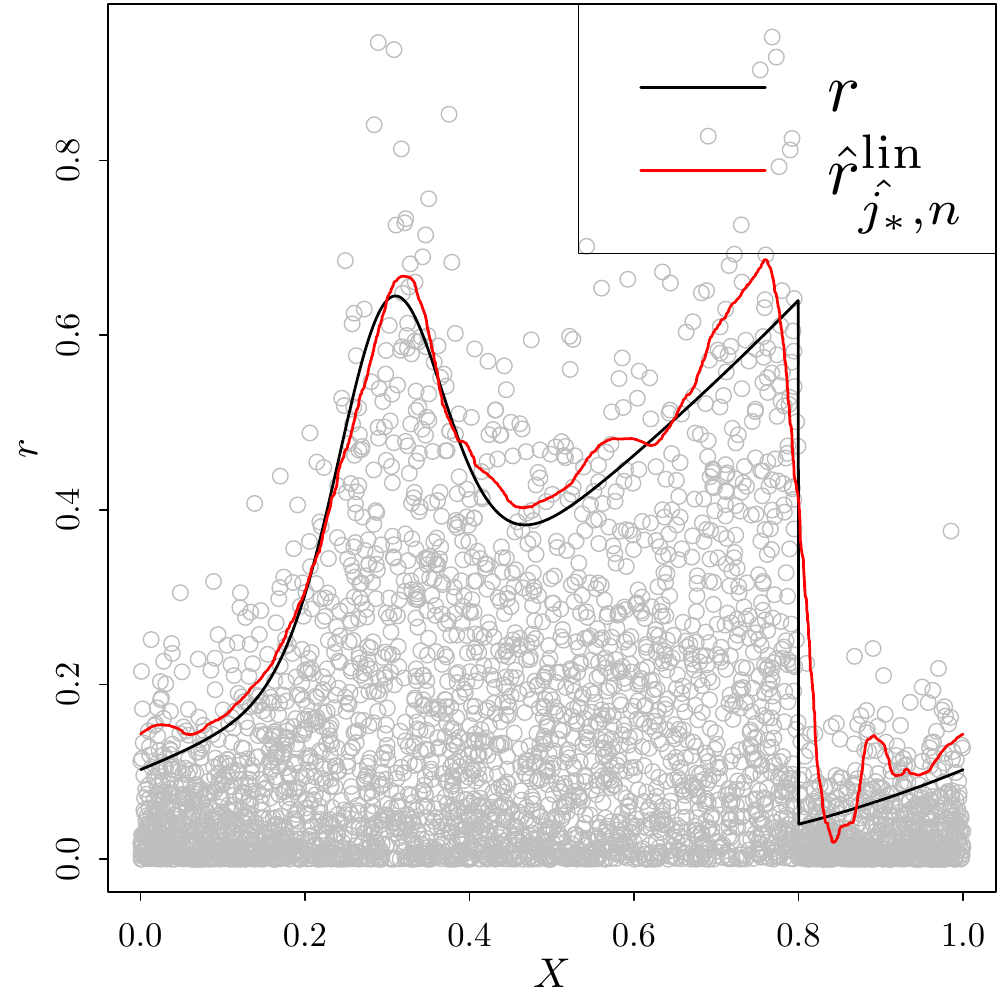}}~
\subfigure[]{
\includegraphics[width=0.31\textwidth]{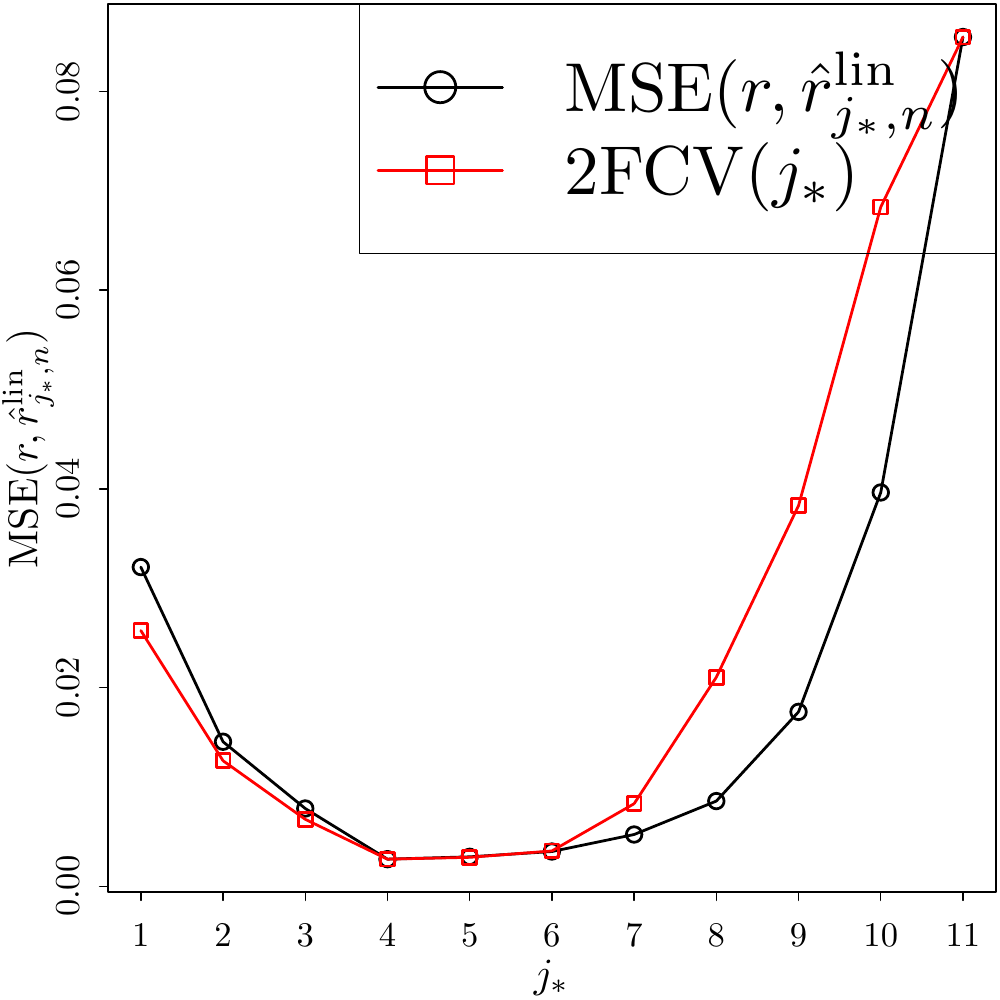}}~
\subfigure[]{
\includegraphics[width=0.31\textwidth]{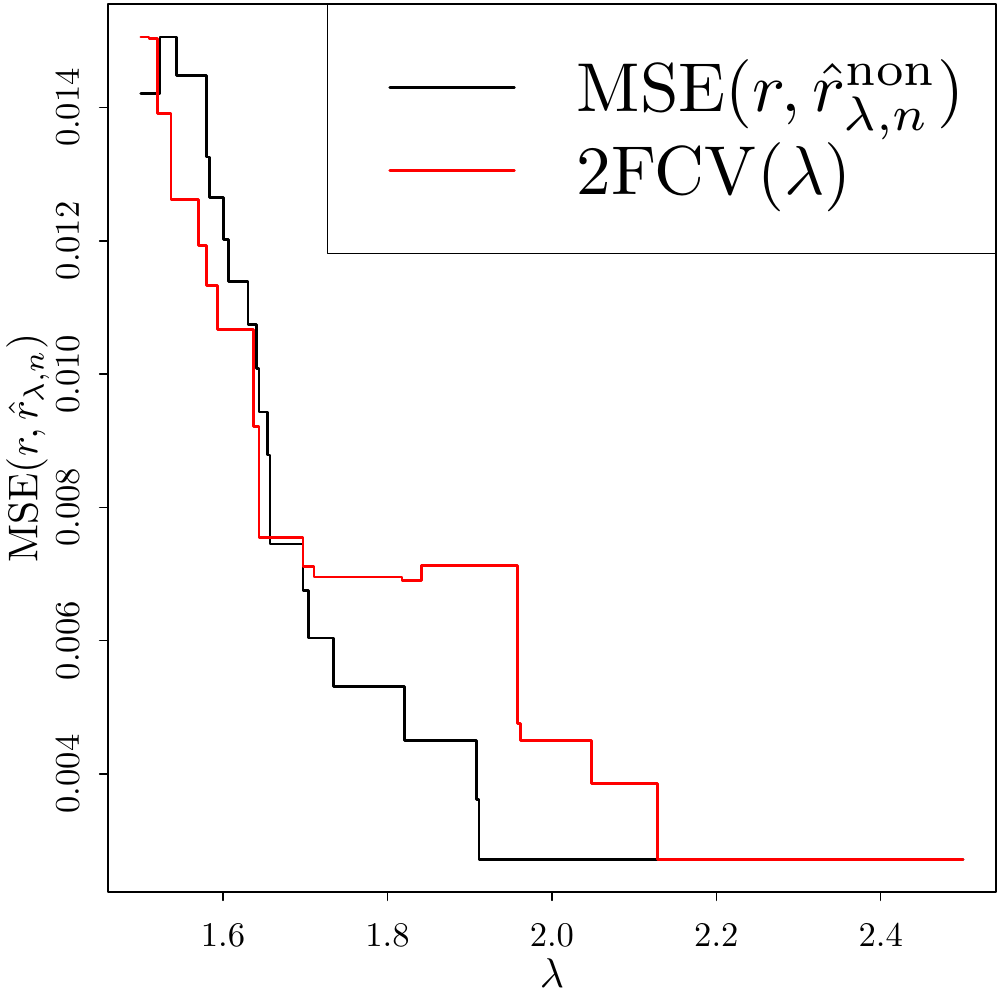}}
\subfigure[]{
\includegraphics[width=0.31\textwidth]{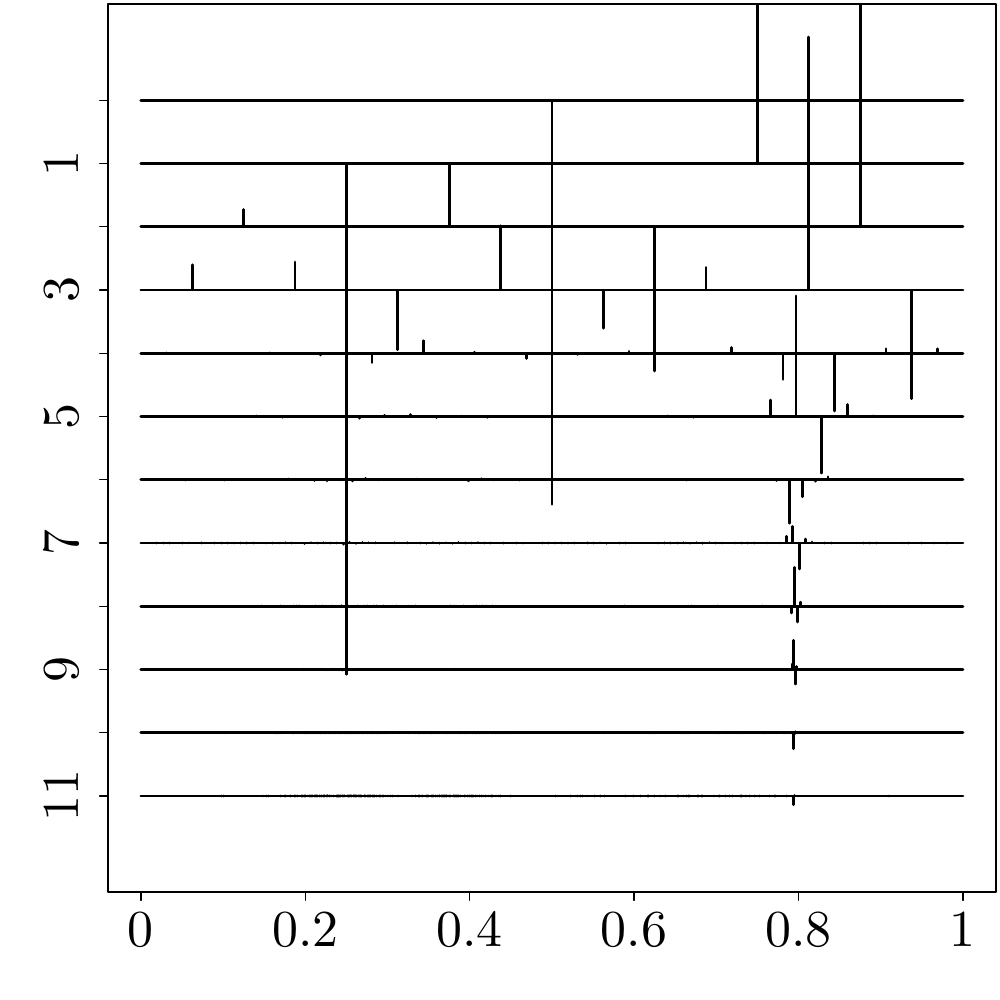}}~
\subfigure[]{
\includegraphics[width=0.31\textwidth]{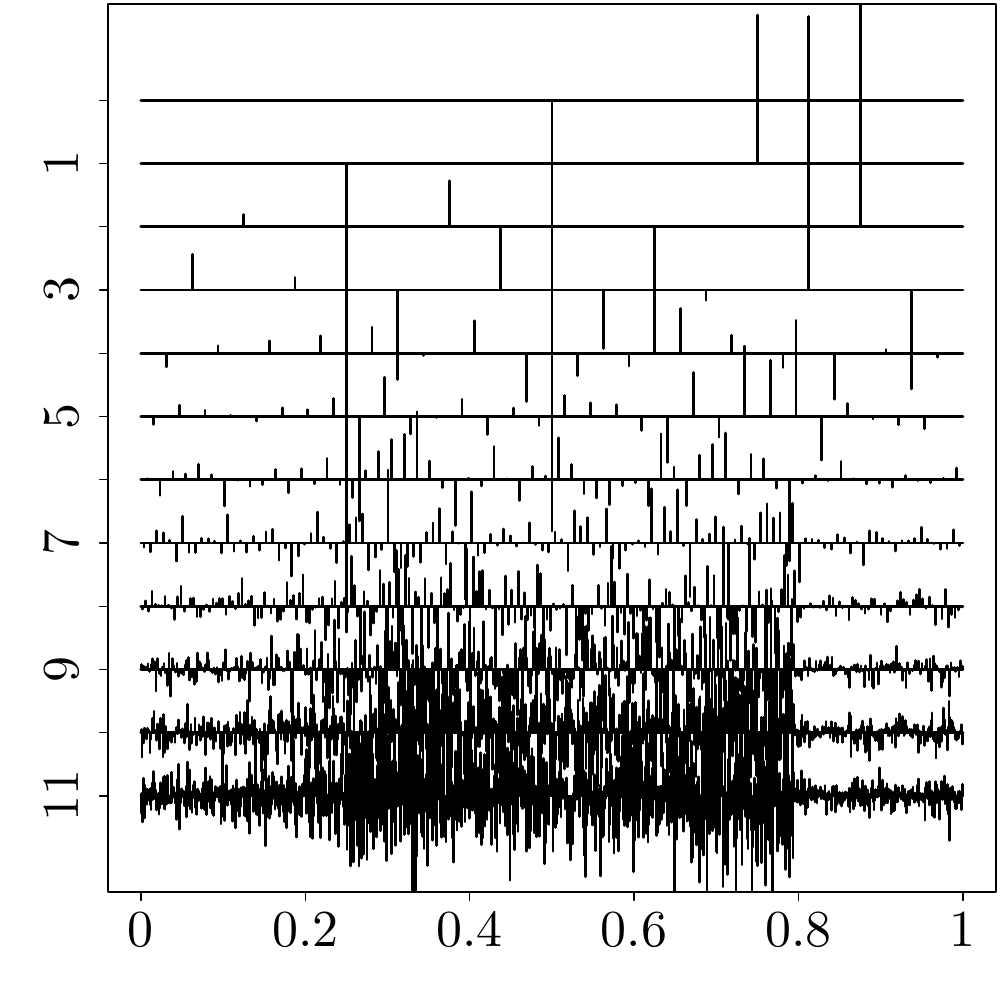}}~
\subfigure[]{
\includegraphics[width=0.31\textwidth]{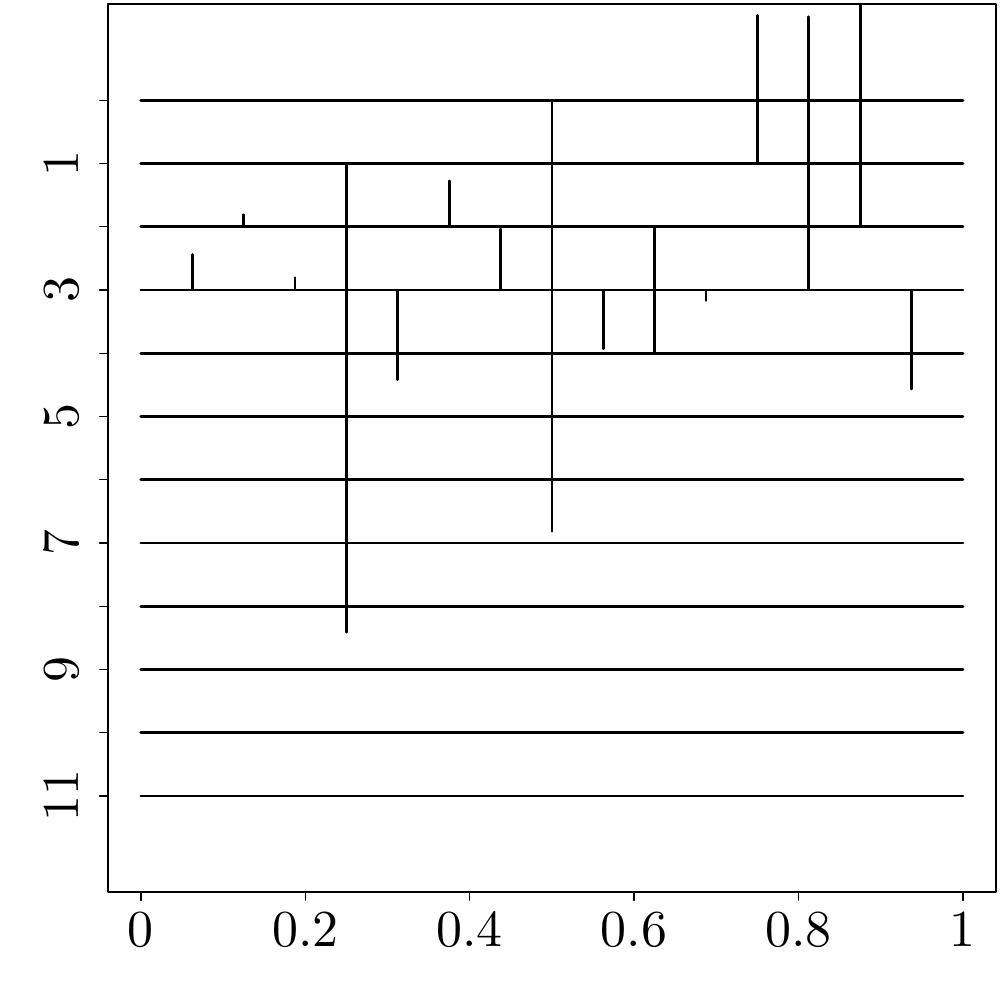}}
\caption{Typical estimation from a single simulation for $n=4096$ and $\sigma^2=0.01$. Noisy observations $(X,Y^2)$ (grey circles), true function (black line) and $\hat{r}_{\hat{j}_*,n}^{\mathrm{lin}}$ (a). Graphs of the MSE (black line) against $j_*$ and (rescaled) 2FCV($\cdot$) criterion (red line) for the linear (b) and nonlinear (c) cases respectively. Original wavelet coefficients (d). Noisy wavelet coefficients (e). Estimated wavelet coefficients (f).}
\label{fig:single}
\end{figure}

A typical example of estimation for the Blip function, with $n=4096$ is given in Figure~\ref{fig:single}. It can be seen that the minimum of $\mathrm{2FCV}(j_*)$ criteria coincides with that of unknown risk (\emph{i.e.}, $\hat{j}_*=4$) and therefore provides the best possible linear estimator for the collection under consideration (\emph{i.e.}, $\mathrm{MSE}(r,\hat{r}_{\hat j_*,n}^{\mathrm{lin}})=0.0027$). We have not included the results of the non-linear estimator in Figure~\ref{fig:single}. Indeed, in this case, the value of the threshold obtained by minimizing the cross validated criterion leads to the elimination of all the thresholded coefficients (\emph{i.e.}, going from $\hat{j}_*$ to $j_1$)  and therefore leads to the same estimate and the same risk as the linear estimator. We can see (Figure~\ref{fig:single}(c)) that the unknown risk behaves in the same way here. This is partly because the amplitude of the coefficients at the fine scales is so large and variable from one scale to another that it is not possible to obtain an overall optimal threshold value that makes it possible to maintain certain important coefficients and that, on the contrary, keeps coefficients associated with noise, with this specific thresholding policy (\emph{i.e.}, a `keep' or `kill' rule). In particular the important coefficients located on scales larger than $\hat{j}_*$ (especially those encoding the discontinuity of $r$) are too small in amplitude to be maintained by a global threshold. 

\begin{figure}[t]
\centering
\subfigure[\textit{Parabolas}]{\includegraphics[width=0.32\textwidth]{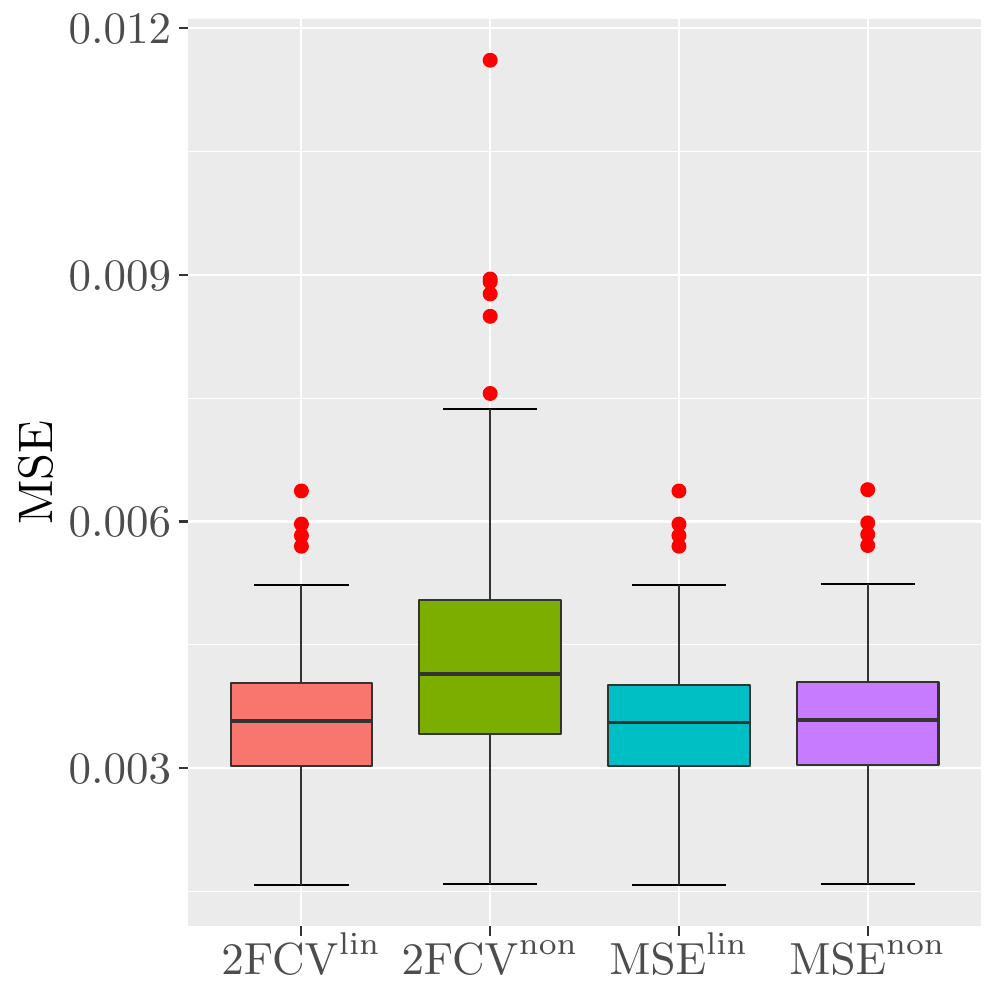}}~
\subfigure[\textit{Ramp}]{\includegraphics[width=0.32\textwidth]{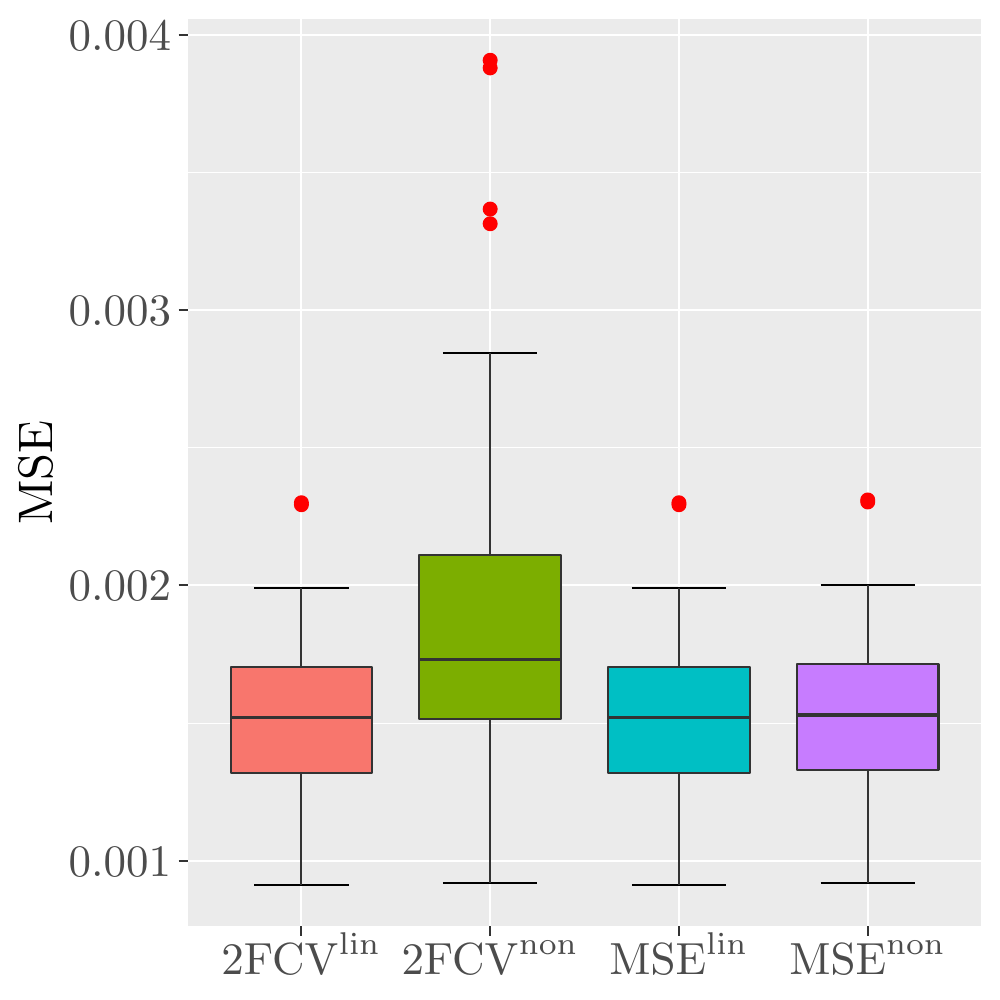}}~
\subfigure[\textit{Blip}]{\includegraphics[width=0.32\textwidth]{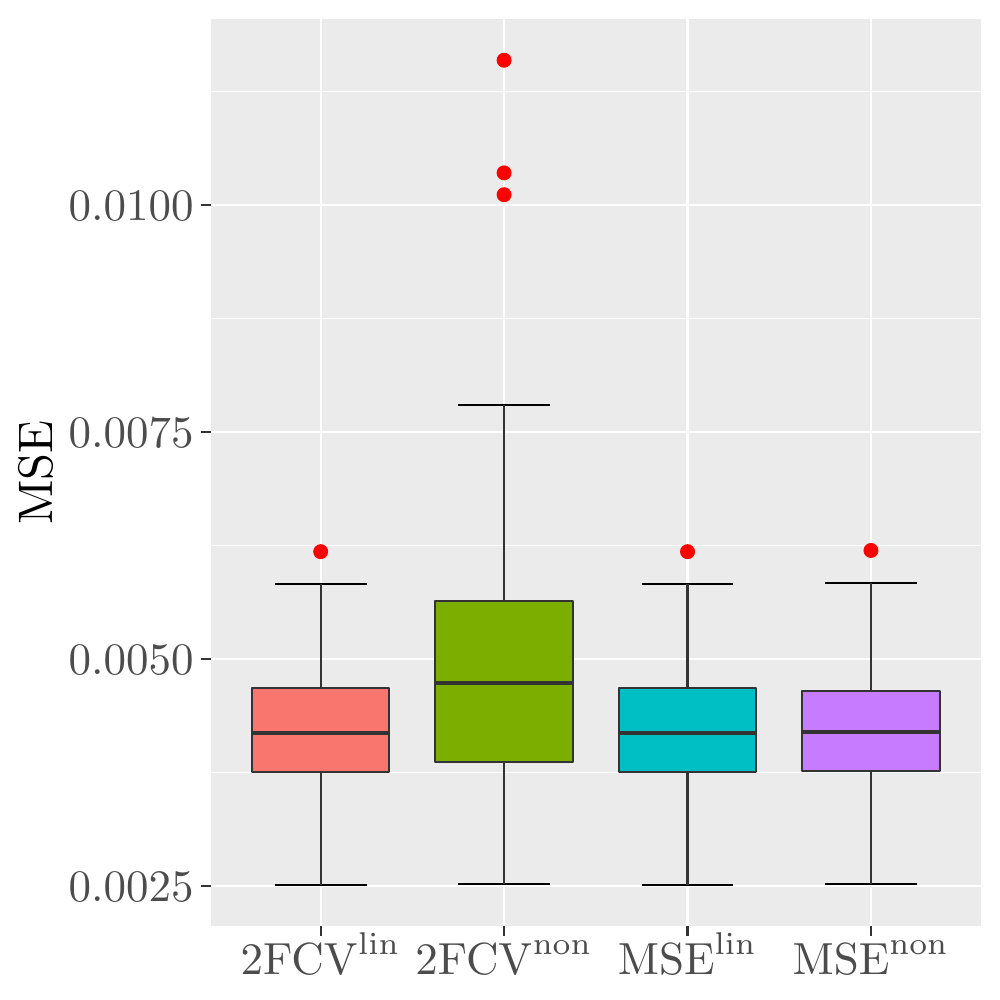}}
\caption{Box plots of the MSE for the three-test functions with $n=4096$ and $\sigma^2=0.01$.}
\label{fig:box}
\end{figure}

In order to determine whether this phenomenon observed for a single function and a single realization is confirmed in a more general context, we compare the performance in terms of MSE (computed on the functions after reconstruction) for both estimators and for the three-test functions. For each function, a sample of $N=100$ is generated and we compare the average behavior of the MSE for parameters selected with the oracle obtained by minimizing the MSE using the original signal $r$ (denoted by $\mathrm{MSE}^{\mathrm{lin}}$ and $\mathrm{MSE}^{\mathrm{non}}$ respectively), the linear $\mathrm{2FCV}^{\mathrm{lin}}$ strategy and the non-linear $\mathrm{2FCV}^{\mathrm{non}}$ (\emph{i.e.}, calculated from $\hat{j}_*$ and the threshold that minimizes the $\mathrm{2FCV}(\lambda)$) strategy. Figure~\ref{fig:box} presents the results in the form of boxplots, one for each function. On the one hand, for all three functions, we can see that the performance of $\mathrm{2FCV}^{\mathrm{lin}}$  is at the $\mathrm{MSE}^{\mathrm{lin}}$ level. This procedure therefore provides a remarkable surrogate of the unknown risk. On the other hand, the non-linear $\mathrm{2FCV}^{\mathrm{non}}$ oracle is similar to $\mathrm{MSE}^{\mathrm{lin}}$, which means that the optimal threshold here leads systematically to the suppression of all threshold coefficients --- which corresponds to the selection of values of the threshold parameter which is greater than the largest noisy wavelet coefficient in absolute value. Finally, the variability of the $\mathrm{2FCV}^{\mathrm{non}}$ is high as a result of selected threshold values that are sometimes too small, resulting in the conservation of unnecessary coefficients in the reconstruction. This is because the curves associated with non-linear criteria do not generally allow a single global minimum, but the minimum is reached in the form of a plateau (see Figure~\ref{fig:single}(c)). In practice, when the minimum is reached on such a plateau, the first element that constitutes it is selected first. This has no influence on $\mathrm{MSE}^{\mathrm{non}}$ but generates this variability of $\mathrm{2FCV}^{\mathrm{non}}$, \emph{i.e.}, when the abscissa of the first point constituting the plateau associated with $\mathrm{2FCV}^{\mathrm{non}}$ is lower than that of $\mathrm{MSE}^{\mathrm{non}}$) Note that to overcome this problem, in the presence of a plateau, we could for example select a threshold value in the middle of it. We have not done so here to emphasize the fact that a cross validation strategy of the global threshold seems ineffective in this setting. It should also be noted that in our simulations, this finding is also verified for other noise levels or sample sizes (the results are generally very similar, so we give only another example by considering a lower number of samples, $n=2048$ and a lower additive noise level $\sigma^2=0.025$).

\begin{figure}[t]
\centering
\subfigure[\textit{Parabolas}]{\includegraphics[width=0.32\textwidth]{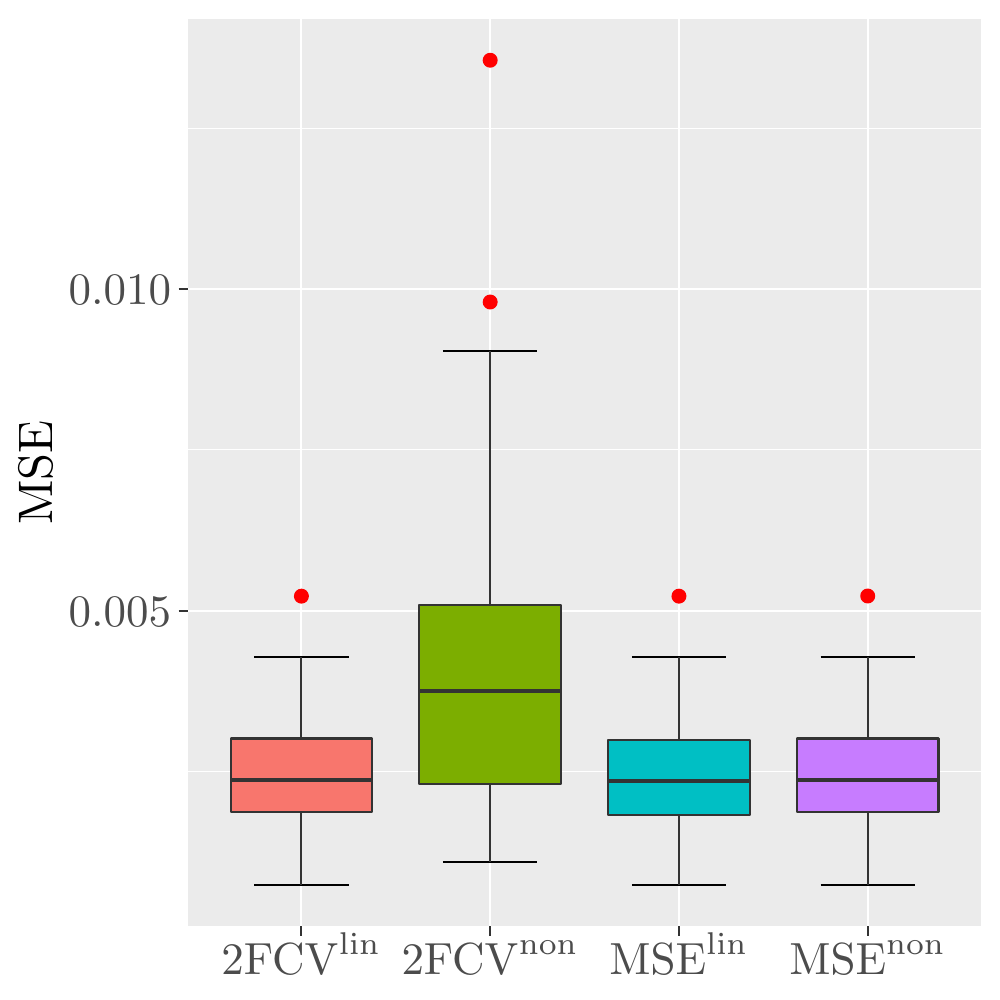}}~
\subfigure[\textit{Ramp}]{\includegraphics[width=0.32\textwidth]{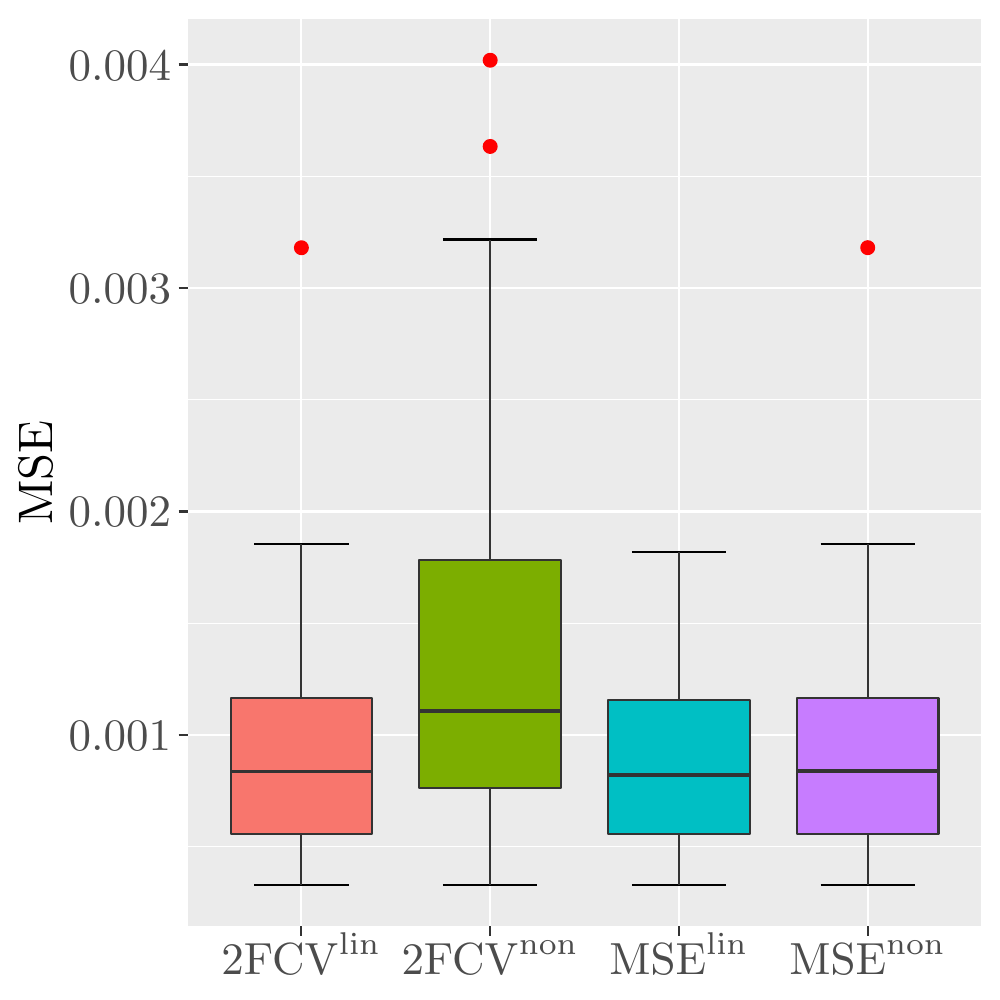}}~
\subfigure[\textit{Blip}]{\includegraphics[width=0.32\textwidth]{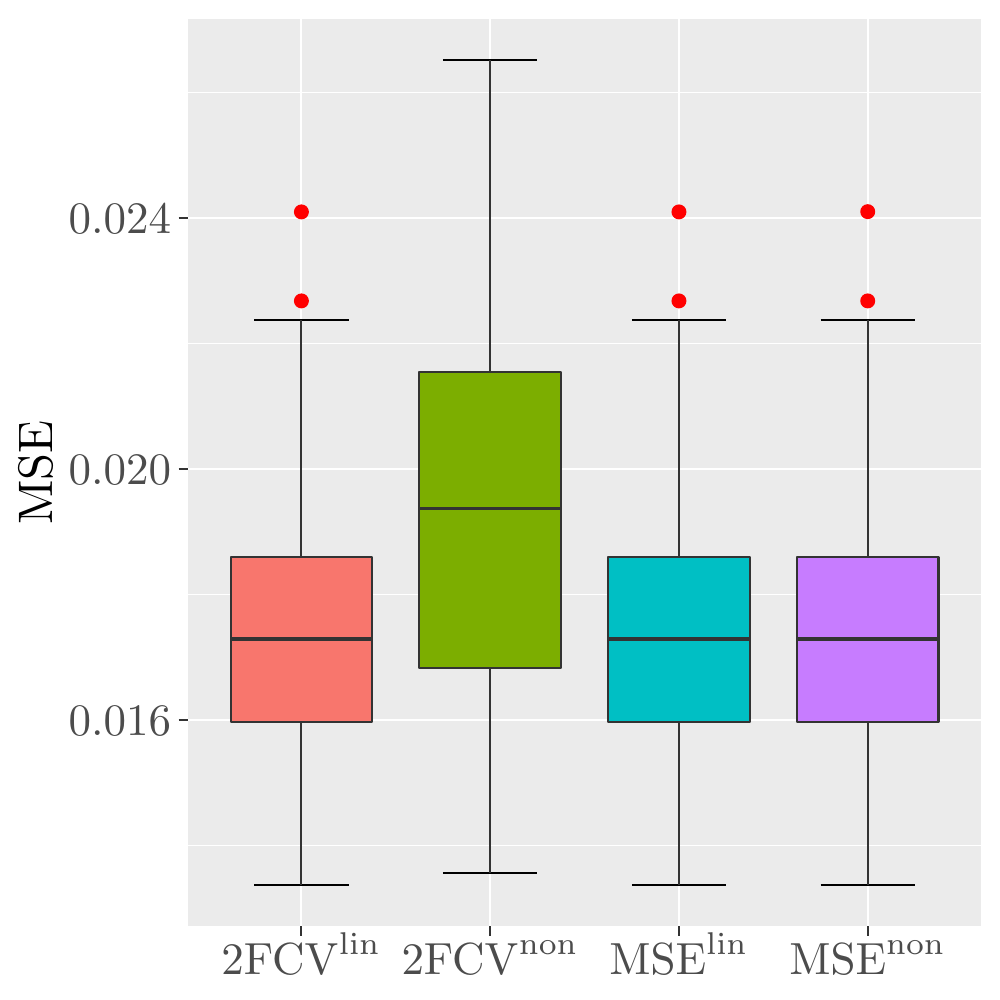}}
\caption{Box plots of the MSE for the three-test functions with $n=2048$ and $\sigma^2=0.025$.}
\label{fig:box2}
\end{figure}

In conclusion, the linear approach seems more appropriate than the non-linear approach in the context of the simulations considered in this study. One way to fully benefit from the non-linear approach would be to consider an optimal threshold selection strategy on a scale by scale basis. The selection procedure used here, based on an interpolation performed in the original domain, does not facilitate this extension. For this purpose it would be necessary, for example, to define an interpolated version of the cross validation method in the wavelet coefficients domain. 

\section{Auxiliary results and proof of the main result}\label{proof}
\subsection{Auxiliary results}\label{lemma}
In this section, we provide some lemmas for the proof of the main Theorem.\par

\begin{lemma}\label{unbiased} Let $j\ge \tau$, $\bk\in \Lambda_j$, $\hat \alpha_{j,\bk}$ be \eqref{ru}. Then, under \ref{hyp:H1} or \ref{hyp:H2}, we have 
\[
\E[\hat{\alpha}_{j,\bk}]={\alpha}_{j,\bk}, \quad \E\left[\frac{1}{n}\sum_{i=1}^{n}Y_i^2\Psi_{j,\bk,u}(\bX_i) -  w_{j,\bk,u} \right]={\beta}_{j,\bk,u}.
\]
\end{lemma}
\noindent{\bf Proof of Lemma \ref{unbiased}.} Using the independence assumptions on the random variables,  \ref{hyp:H1} or \ref{hyp:H2}, observe that
\begin{align*}
\E\left[U_1 V_1 f  (\bX_1)\Phi_{j,\bk}(\bX_1) \right]&=
\begin{cases}
 \E[U_1]\E[V_1]\E\left[f  (\bX_1)\Phi_{j,\bk}(\bX_1) \right] & \text{under \ref{hyp:A5}},\\
 \E[U_1]\E\left[V_1f  (\bX_1)\Phi_{j,\bk}(\bX_1) \right]      & \text{under \ref{hyp:A6}},
\end{cases}\\
&=0,
\end{align*}
and 
\begin{align*}
v_{j,\bk}&=
\begin{cases}
\E\left[V_1^2\right] 2^{-jd/2}=\E\left[V_1^2\right] \int_{[0,1]^d} \Phi_{j,\bk}(\bx)d\bx=\E\left[V_1^2\right]  \E\left[\Phi_{j,\bk}(\bX_1)\right]& \text{under \ref{hyp:A5}},\\
\displaystyle \int_{[0,1]^d}g^2(\bx) \Phi_{j,\bk}(\bx) d\bx& \text{under \ref{hyp:A6}},
\end{cases}\\
&= \E\left[V_1^2\Phi_{j,\bk}(\bX_1)\right].
\end{align*}
Therefore
\begin{align*}
\E[\hat{\alpha}_{j,\bk}]&=\E\left[\frac{1}{n}\sum_{i=1}^{n}Y_i^2\Phi_{j,\bk}(\bX_i) - v_{j,\bk}\right]= \E\left[Y_1^2\Phi_{j,\bk}(\bX_1)\right] -  v_{j,\bk}\\
& = \E\left[U_1^2r (\bX_1)\Phi_{j,\bk}(\bX_1)\right]+2\E\left[U_1 V_1 f  (\bX_1)\Phi_{j,\bk}(\bX_1) \right]+  \E\left[V_1^2\Phi_{j,\bk}(\bX_1)\right]-v_{j,\bk}\\
& =\E\left[U_1^2 \right] \E\left[r (\bX_1)\Phi_{j,\bk}(\bX_1)\right]= \int_{[0,1]^d} r(\bx)\Phi_{j,\bk}(\bx)d\bx={\alpha}_{j,\bk}.
\end{align*}
Using similar mathematical arguments, since $\int_{[0,1]^d} \Psi_{j,\bk,u}(\bx)d\bx=0$, we have 
\begin{align*}
w_{j,\bk,u}&=
\begin{cases}
 \displaystyle0= \E\left[V_1^2\right]\int_{[0,1]^d} \Psi_{j,\bk,u}(\bx)d\bx=\E\left[V_1^2\right]\E\left[\Psi_{j,\bk,u}(\bX_1)\right]& \text{under \ref{hyp:A5}},\\
\displaystyle\int_{[0,1]^d}g^2(\bx) \Psi_{j,\bk,u}(\bx) d\bx& \text{under \ref{hyp:A6}},
\end{cases}\\
&= \E\left[V_1^2\Psi_{j,\bk,u}(\bX_1)\right].
\end{align*}
We prove the second equality. The proof of Lemma \ref{unbiased} is complete. $\hfill \Box$

\begin{lemma}\label{var2} Let $j\ge \tau$ such that $2^j\le n$, $\bk\in \Lambda_j$, $\hat \alpha_{j,\bk}$ and $\hat \beta_{j,\bk,u}$ be \eqref{ru} and \eqref{ru0} respectively. Then, under \ref{hyp:H1} or \ref{hyp:H2}, 
\begin{equation*}
\E\left[ (\hat{\alpha}_{j,\bk}-\alpha_{j,\bk} )^2\right]\lesssim    \frac{1}{n}, \quad  \E\left[ (\hat{\beta}_{j,\bk,u}-\beta_{j,\bk,u} )^2\right]\lesssim   \frac{\ln n}{n}.
\end{equation*}
\end{lemma}
\noindent{\bf Proof of Lemma \ref{var2}.} Owing to Lemma \ref{unbiased} we have $\E[\hat{\alpha}_{j,\bk}]={\alpha}_{j,\bk}$. Therefore
\begin{align}\label{irre}
\E\left[ (\hat{\alpha}_{j,\bk}-\alpha_{j,\bk} )^2\right] &= \Var\left[ \hat{\alpha}_{j,\bk}\right]=\Var\left[ \frac{1}{n}\sum_{i=1}^{n}Y_i^2\Phi_{j,\bk}(\bX_i) - v_{j,\bk}\right]\nonumber\\
& = \Var\left[ \frac{1}{n}\sum_{i=1}^{n}Y_i^2\Phi_{j,\bk}(\bX_i) \right]=\frac{1}{n} \Var\left[Y_1^2\Phi_{j,\bk}(\bX_1) \right]\le \frac{1}{n} \E\left[Y_1^4\Phi^2_{j,\bk}(\bX_1) \right]\nonumber\\
& \lesssim \frac{1}{n} \left(\E\left[U_1^4f^4(\bX_1)\Phi^2_{j,\bk}(\bX_1) \right] + \E\left[V_1^4\Phi^2_{j,\bk}(\bX_1) \right] \right)\nonumber\\
& = \frac{1}{n} \left(\E\left[U_1^4\right] \E\left[f^4(\bX_1)\Phi^2_{j,\bk}(\bX_1) \right] + \E\left[V_1^4\Phi^2_{j,\bk}(\bX_1) \right] \right).
\end{align}
By \ref{hyp:A1} and $\E\left[\Phi^2_{j,\bk}(\bX_1) \right]=\int_{[0,1]^d} \left(\Phi_{j,\bk}(\bx) \right)^2d\bx=1$, we have $\E\left[f^4(\bX_1)\Phi^2_{j,\bk}(\bX_1) \right]\lesssim 1$. On the other hand, we have 
\begin{equation*}
\E\left[V_1^4\Phi^2_{j,\bk}(\bX_1) \right]=
\begin{cases}
 \displaystyle \E\left[V_1^4\right] \E\left[\Phi^2_{j,\bk}(\bX_1) \right]\lesssim 1& \text{under \ref{hyp:A5}},\\
\displaystyle\int_{[0,1]^d}g^4(\bx) \Phi^2_{j,\bk}(\bx) d\bx\lesssim \int_{[0,1]^d}\Phi^2_{j,\bk}(\bx) d\bx=1& \text{under \ref{hyp:A6}}.
\end{cases}
\end{equation*}
Thus all the terms in the brackets of \eqref{irre} are bounded from above. The first inequality in Lemma \ref{var2} is proved. 

Now, by the definition of $\hat{\beta}_{j,\bk,u}$, taking $K_i:=Y_i^2\Psi_{j,\bk,u}(\bX_i) - w_{j,\bk,u}$ and $D_{i}:=K_i\bbone_{\{|K_i|\leq\rho_n\}}-\E\left[ K_i\bbone_{\{|K_i|\leq\rho_n\}}\right]$ for the sake of simplicity, the second equation in Lemma \ref{unbiased} yields 
\begin{equation*}
\hat{\beta}_{j,\bk,u}-\beta_{j,\bk,u}=\frac{1}{n}\sum_{i=1}^{n}D_{i}-\E\left[K_1\bbone_{\{|K_1|>\rho_n\}}\right].
\end{equation*}
Hence, using $\E\left[ \left((1/n) \sum\limits_{i=1}^{n}D_{i}\right)^2\right]= ({1}/{n^2})\Var\left[\sum\limits_{i=1}^{n}D_{i}\right]=\Var[D_1]/n\le \E\left[D_1^2\right]/n$, we have 
\begin{align*}
\E\left[(\hat{\beta}_{j,\bk,u}-\beta_{j,\bk,u})^{2}\right]& \lesssim  \E\left[ \left(\frac{1}{n} \sum\limits_{i=1}^{n}D_{i}\right)^2\right]+\left(\E\left[|K_1|\bbone_{\{|K_1|>\rho_n\}}\right]\right)^{2}.\\
& \lesssim  \frac{1}{n}\E \left[D^{2}_1\right]+\left(\E\left[|K_1|\bbone_{\{|K_1|>\rho_n\}}\right]\right)^{2}.
\end{align*}

Proceeding as for the proof of the first inequality, using the assumptions \ref{hyp:H1} or \ref{hyp:H2}, note that $\E\left[K^{2}_{1}\right]\lesssim \E\left[Y^{4}_{1}\Psi^{2}_{j,\bk,u}(\bX_1)\right] + w^2_{j,\bk,u}\lesssim1$ and 
\begin{equation}\label{lm24}
\left(\E\left[K_1\bbone_{\{|K_1|>\rho_n\}}\right]\right)^{2}\lesssim\left(\E \left[K_1^2 \right]/\rho_n\right)^{2}\lesssim \frac{\ln n}{n}.
\end{equation}
Therefore,
\[
\E\left[(\hat{\beta}_{j,\bk,u}-\beta_{j,\bk,u})^{2}\right]\lesssim \frac{1}{n}+\frac{\ln n}{n}\lesssim\frac{\ln n}{n}.
\]
The second  inequality in Lemma \ref{var2} is proved. This ends the proof of Lemma \ref{var2}. $\hfill \Box$

\begin{lemma}\label{var3} Let $j\ge \tau$ such that $2^{jd}\lesssim n/\ln n$, $\bk\in \Lambda_j$, $\hat \beta_{j,\bk,u}$ be \eqref{ru0}. Then there exists a constant $\kappa>1$ such that
\begin{equation*}
\mathbb{P}(|\hat{\beta}_{j,\bk,u}-\beta_{j,\bk,u}|\geq\kappa t_{n})\lesssim n^{-4}.
\end{equation*}
\end{lemma}

\noindent{\bf Proof of Lemma \ref{var3}.} By the definition of $\hat{\beta}_{j,\bk,u}$, taking $K_i:=Y_i^2\Psi_{j,\bk,u}(\bX_i) - w_{j,\bk,u}$ and $D_{i}:=K_i\bbone_{\{|K_i|\leq\rho_n\}}-E\left[ K_i\bbone_{\{|K_i|\leq\rho_n\}}\right]$ for the sake of simplicity, the second equation in Lemma \ref{unbiased} yields 
\begin{equation*}
|\hat{\beta}_{j,\bk,u}-\beta_{j,\bk,u}|\lesssim\frac{1}{n}\left|\sum_{i=1}^{n}D_{i}\right|+E\left[ |K_1|\bbone_{\{|K_1|>\rho_n\}}\right].
\end{equation*}
Using \eqref{lm24}, there exists $c>0$ such that $\E\left[ |K_1|\bbone_{\{|K_1|>\rho_n\}}\right]\le c \sqrt{\ln n/n}$. Then
\begin{equation*}
\{|\hat{\beta}_{j,\bk,u}-\beta_{j,\bk,u}|\geq\kappa t_{n}\}\subseteq\left\lbrace \frac{1}{n} \left|\sum_{i=1}^{n}D_{i}\right|\geq(\kappa-c) t_{n}\right\rbrace.
\end{equation*}
Note that $\E[D_{i}]=0$ thanks to Lemma \ref{unbiased}. According to the proof of Lemma \ref{var2}, $\E[D_{i}^{2}]:=\delta^{2}\lesssim1$. This with $|D_{i}|\lesssim \sqrt{n/\ln n}$ and Bernstein inequality shows
\begin{align*}
\mathbb{P}\left( \frac{1}{n}\left|\sum_{i=1}^{n}D_{i}\right|\geq(\kappa-c) t_{n}\right)&\lesssim\exp{\left\lbrace-\frac{n(\kappa-c)^{2}t^{2}_{n}}{2(\delta^{2}+(\kappa-c) t_{n}\rho_n/3)}\right\rbrace}\notag\\
&\lesssim\exp{\left\lbrace-\ln n\frac{(\kappa-c)^{2}}{2(\delta^{2}+(\kappa-c)/3)}\right\rbrace}\notag\\
&\lesssim n^{-\frac{(\kappa-c)^{2}}{2(\delta^{2}+(\kappa-c)/3)}}.
\end{align*}
Then one choose large enough $\kappa$ such that
\begin{equation}
\mathbb{P}(|\hat{\beta}_{j,\bk,u}-\beta_{j,\bk,u}|\geq\kappa t_{n})
\lesssim n^{-\frac{(\kappa-1)^{2}}{2(\delta^{2}+(\kappa-1)/3)}}\lesssim n^{-4}.\notag
\end{equation}
This is the desired conclusion.$\hfill \Box$

\subsection{Proof of the main result}
This section is devoted to the proof of Theorem \ref{theo}. We prove \eqref{eq11a} and \eqref{eq11b} in turn. 

\noindent \textbf{Proof of \eqref{eq11a}}  ~Note that
\begin{equation}\label{eqq1}
\E\left[\int_{[0,1]^{d}}\big|\hat{r}^{\mathrm{lin}}_{n}(\bx)-r(\bx)\big|^{2}d\bx\right]= \E\left[\big\|\hat{r}^{\mathrm{lin}}_{n}-P_{j_{*}}r\big\|_{2}^{2}\right]+\big\|P_{j_{*}}r-r\big\|_{2}^{2}.
\end{equation}
It is easy to see that
\[
\E\left[\left\|\hat{r}^{\mathrm{lin}}_{n}-P_{j_{*}}r\right\|^{2}_{2}\right]=\E\left[\left\|\sum\limits_{\bk\in\Lambda_{j_{*}}}(\hat{\alpha}_{j_{*},\bk}-\alpha_{j_{*},\bk})\Phi_{j_{*},\bk}\right\|^{2}_{2}\right] =\sum\limits_{\bk\in\Lambda_{j_{*}}}
\E\left[\Big|\hat{\alpha}_{j_{*},\bk}-\alpha_{j_{*},\bk}\Big|^{2}\right].
\]
According to Lemma \ref{var2}, $|\Lambda_{j_{*}}|\thicksim2^{j_{*}d}$ and $2^{j_{*}}\thicksim n^{\frac{1}{2s'+d}}$,
\begin{equation}\label{eqq2}
\E\left[\left\|\hat{r}^{\mathrm{lin}}_{n}-P_{j_{*}}r\right\|_{2}^{2}\right]\lesssim \frac{2^{j_{*}d}}{n}\thicksim n^{-\frac{2s'}{2s'+d}}.
\end{equation}
When $p\geq2$, $s'=s$. By H\"{o}lder inequality and $r\in B_{p,q}^{s}([0,1]^{d})$,
\begin{equation*}
\|P_{j_{*}}r-r\|_{2}^{2}\lesssim\|P_{j_{*}}r-r\|_{p}^{2}\lesssim2^{-2j_{*}s}\thicksim n^{-\frac{2s}{2s+d}}.
\end{equation*}
When $1\leq p<2$ and $s>d/p$, $B_{p,q}^{s}([0,1]^{d})\subseteq B_{2,\infty}^{s'}([0,1]^{d})$ 
\begin{equation*}
\|P_{j_{*}}r-r\|_{2}^{2}\lesssim\sum\limits_{j=j_{*}}^{\infty}2^{-2js'}\lesssim2^{-2j_{*}s'}\thicksim n^{-\frac{2s'}{2s'+d}}.
\end{equation*}
Therefore, in both cases, 
\begin{equation}\label{eqq5}
\|P_{j_{*}}r-r\|_{2}^{2}\lesssim n^{-\frac{2s'}{2s'+d}}.
\end{equation}
By \eqref{eqq1}, \eqref{eqq2} and \eqref{eqq5},
\begin{equation*}
\E\left[\int_{[0,1]^{d}}\big|\hat{r}^{\mathrm{lin}}_{n}(\bx)-r(\bx)\big|^{2}d\bx\right]\lesssim n^{-\frac{2s'}{2s'+d}}.\notag
\end{equation*}

\noindent \textbf{Proof of \eqref{eq11b}}  We now follow the lines of \cite[Theorem 2]{delyon} with adaptation to our statistical setting, by using the definitions of our estimators and the auxiliary results of Section \ref{lemma}. By the definitions of $\hat{r}^{\mathrm{lin}}_{n}$ and $\hat{r}^{\mathrm{non}}_{n}$, we have 
\begin{align*}
\hat{r}^{\mathrm{non}}_{n}(\bx)-r(\bx)=&\Big(\hat{r}^{\mathrm{lin}}_{n}(\bx)-P_{j_{*}}r(\bx)\Big)-\Big(r(\bx)-P_{j_{1}+1}r(\bx)\Big)
\\
& + \sum\limits_{j=j_{*}}^{j_{1}} \sum\limits_{u=1}^{2^{d}-1}\sum\limits_{\bk\in\Lambda_j}\Big(\hat{\beta}_{j,\bk,u}\bbone_{\{|\hat{\beta}_{j,\bk,u}|\geq\kappa t_{n}\}}-\beta_{j,\bk,u}\Big)\Psi_{j,\bk,u}(\bx).
\end{align*}
  Hence,
\begin{equation*}
\E\left[\int_{[0,1]^{d}}\Big|\hat{r}^{\mathrm{non}}_{n}(\bx)-r(\bx)\Big|^{2}d\bx\right]\lesssim T_{1}+T_{2}+Q,
\end{equation*}
where $T_{1}:=\E\left[\Big\|\hat{r}^{\mathrm{lin}}_{n}-P_{j_{*}}r\Big\|^{2}_{2}\right],~~T_{2}:=\Big\|r-P_{j_{1}+1}r\Big\|^{2}_{2}$ and
\[
Q:=\E\left[\left\|\sum\limits_{j=j_{*}}^{j_{1}} \sum\limits_{u=1}^{2^{d}-1}\sum\limits_{\bk\in\Lambda_j}\left(\hat{\beta}_{j,\bk,u}\bbone_{\{|\hat{\beta}_{j,\bk,u}|\geq\kappa t_{n}\}}-\beta_{j,\bk,u}\right)\Psi_{j,\bk,u}\right\|^{2}_{2}\right].
\]
According to \eqref{eqq2} and $2^{j_{*}}\sim n^{\frac{1}{2m+d}}~(m>s)$,
\[
 T_{1}=\E\left[\Big\|\hat{r}^{\mathrm{lin}}_{n}-P_{j_{*}}r\Big\|_{2}^{2}\right]\lesssim \frac{2^{j_{*}d}}{n}\thicksim n^{-\frac{2m}{2m+d}}<n^{-\frac{2s}{2s+d}}.
 \]
When $p\geq2$, by the same arguments as \eqref{eqq5} shows
$T_{2}=\Big\|r-P_{j_{1}+1}r\Big\|^{2}_{2}\lesssim2^{-2j_{1}s}.$ This with $2^{j_{1}}\sim(n/\ln n)^{\frac{1}{d}}$ leads to
\begin{equation*}
T_{2}\lesssim2^{-2j_{1}s}\thicksim\left(\frac{\ln n}{n}\right)^{\frac{2s}{d}}\leq (\ln n)n^{-\frac{2s}{2s+d}}.
\end{equation*}
On the other hand, $B_{p,q}^{s}([0,1]^{d})\subseteq B_{2,\infty}^{s+d/2-d/p}([0,1]^{d})$ when $1\leq p<2$ and $s>d/p$. Then
\begin{equation}
T_{2}\lesssim2^{-2j_{1}(s+\frac{d}{2}-\frac{d}{p})}\thicksim\big(\frac{\ln n}{n}\big)^{\frac{2(s+\frac{d}{2}-\frac{d}{p})}{d}}\leq (\ln n)n^{-\frac{2s}{2s+d}}.\notag
\end{equation}
Hence,
\[
T_{2}\lesssim(\ln n)n^{-\frac{2s}{2s+d}},
\]
for each $1\leq p<+\infty$.

The main work for the proof of \eqref{eq11b} is to show
\begin{equation*}
Q=\E\left[\left\|\sum\limits_{j=j_{*}}^{j_{1}} \sum\limits_{u=1}^{2^{d}-1}\sum\limits_{\bk\in\Lambda_j}\left(\hat{\beta}_{j,\bk,u}\bbone_{\{|\hat{\beta}_{j,\bk,u}|\geq\kappa t_{n}\}}-\beta_{j,\bk,u}\right)\Psi_{j,\bk,u}\right\|^{2}_{2}\right]\lesssim(\ln n){n}^{-\frac{2s}{2s+d}}.
\end{equation*}
Note that
\begin{equation}\label{q2}
Q=\sum\limits_{j=j_{*}}^{j_{1}}\sum\limits_{u=1}^{2^{d}-1}\sum\limits_{\bk\in\Lambda_{j}}\E\left[\left|\hat{\beta}_{j,\bk,u}\bbone_{\{|\hat{\beta}_{j,\bk,u}|\geq\kappa t_{n}\}}-\beta_{j,\bk,u}\right|^{2}\right]\lesssim Q_{1}+Q_{2}+Q_{3},
\end{equation}
where 
\[
Q_{1}=\sum\limits_{j=j_{*}}^{j_{1}}\sum\limits_{u=1}^{2^{d}-1}\sum\limits_{\bk\in\Lambda_{j}}\E\left[\left|\hat{\beta}_{j,\bk,u}-\beta_{j,\bk,u}\right|^{2}\bbone_{\{|\hat{\beta}_{j,\bk,u}-\beta_{j,\bk,u}|>\frac{\kappa t_{n}}{2}\}}\right],
\]
\[
Q_{2}=\sum\limits_{j=j_{*}}^{j_{1}}\sum\limits_{u=1}^{2^{d}-1}\sum\limits_{\bk\in\Lambda_{j}}\E\left[\left|\hat{\beta}_{j,\bk,u}-\beta_{j,\bk,u}\right|^{2}\bbone_{\{|\beta_{j,\bk,u}|\geq\frac{\kappa t_{n}}{2}\}}\right],
\]
\[
Q_{3}=\sum\limits_{j=j_{*}}^{j_{1}}\sum\limits_{u=1}^{2^{d}-1}\sum\limits_{\bk\in\Lambda_{j}}\left|\beta_{j,\bk,u}\right|^{2}\bbone_{\{|\beta_{j,\bk,u}|\leq2\kappa t_{n}\}}.
\]
For $Q_{1}$, one observes that
\[
\E\left[\left|\hat{\beta}_{j,\bk,u}-\beta_{j,\bk,u}\right|^{2}\bbone_{\{|\hat{\beta}_{j,\bk,u}-\beta_{j,\bk,u}|>\frac{\kappa t_{n}}{2}\}}\right]\leq\left(\E\left[\left|\hat{\beta}_{j,\bk,u}-\beta_{j,\bk,u}\right|^{4}\right]\right)^{\frac{1}{2}}\left(\mathbb{P} \left(|\hat{\beta}_{j,\bk,u}-\beta_{j,\bk,u}|>\frac{\kappa t_{n}}{2}\right)\right)^{\frac{1}{2}}
\]
thanks to H\"{o}lder inequality. By Lemma \ref{var2}, Lemma \ref{var3} and $|\hat{\beta}_{j,\bk,u}-\beta_{j,\bk,u}|^{2}\lesssim n/\ln n$,
\begin{equation*}
\E\left[\left|\hat{\beta}_{j,\bk,u}-\beta_{j,\bk,u}\right|^{2}\bbone_{\{|\hat{\beta}_{j,\bk,u}-\beta_{j,\bk,u}|>\frac{\kappa t_{n}}{2}\}}\right]\lesssim \frac{1}{n^{2}}.
\end{equation*}
Then $Q_{1}\lesssim\sum\limits_{j=j_{*}}^{j_{1}}2^{jd}/n^{2}\lesssim 2^{j_{1}d} / n^{2}\lesssim 1/n\leq n^{-\frac{2s}{2s+d}}$, where one uses the choice $2^{j_{1}}\sim (n/\ln n)^{\frac{1}{d}}$. Hence,
\begin{equation}\label{q3}
Q_{1}\leq n^{-\frac{2s}{2s+d}}.
\end{equation}
To estimate $Q_{2}$, one defines
\[
2^{j'}\sim n^{\frac{1}{2s+d}}.
\]
It is easy to see that $2^{j_{*}}\sim n^{\frac{1}{2m+d}}\leq2^{j'}\sim n^{\frac{1}{2s+d}}\leq2^{j_{1}}\sim(n/\ln n)^{\frac{1}{d}}$. Furthermore,
one rewrites
\begin{align*}
Q_{2}&=\left(\sum\limits_{j=j_{*}}^{j'}+\sum\limits_{j=j'+1}^{j_{1}}\right)\left\{\sum\limits_{u=1}^{2^{d}-1}\sum\limits_{\bk\in\Lambda_{j}}\E\left[\left|\hat{\beta}_{j,\bk,u}-\beta_{j,\bk,u}\right|^{2}\bbone_{\{|\beta_{j,\bk,u}|\geq\frac{\kappa t_{n}}{2}\}}\right]\right\}\\
&:=Q_{21}+Q_{22}.
\end{align*}
By Lemma \ref{var2} and $2^{j'}\sim n^{\frac{1}{2s+d}}$,
\begin{align*}
Q_{21}&:=\sum\limits_{j=j_{*}}^{j'}\sum\limits_{u=1}^{2^{d}-1}\sum\limits_{\bk\in\Lambda_{j}}\E\left[\left|\hat{\beta}_{j,\bk,u}-\beta_{j,\bk,u}\right|^{2}\bbone_{\{|\beta_{j,\bk,u}|\geq\frac{\kappa t_{n}}{2}\}}\right]\notag\\
&\lesssim\sum\limits_{j=j_{*}}^{j'}\sum\limits_{u=1}^{2^{d}-1}\sum\limits_{\bk\in\Lambda_{j}}\frac{\ln n}{n}
\lesssim\sum\limits_{j=j_{*}}^{j'}(\ln n)\frac{2^{jd}}{n}\lesssim(\ln n)\frac{2^{j'd}}{n}\sim (\ln n)n^{-\frac{2s}{2s+d}}.
\end{align*}
On the other hand, it follows from Lemma \ref{var2} that
\begin{align*}
Q_{22}&:=\sum\limits_{j=j'+1}^{j_{1}}\sum\limits_{u=1}^{2^{d}-1}\sum\limits_{\bk\in\Lambda_{j}}\E\left[\left|\hat{\beta}_{j,\bk,u}-\beta_{j,\bk,u}\right|^{2}\bbone_{\{|\beta_{j,\bk,u}|\geq\frac{\kappa t_{n}}{2}\}}\right]\\
&\lesssim\sum\limits_{j=j'+1}^{j_{1}}\sum\limits_{u=1}^{2^{d}-1}\sum\limits_{\bk\in\Lambda_{j}}\frac{\ln n}{n}~\bbone_{\{|\beta_{j,\bk,u}|\geq\frac{\kappa t_{n}}{2}\}}.
\end{align*}
When $p\geq2$, since $r\in B_{p,q}^{s}([0,1]^{d})$, Lemma \ref{var2} and $t_{n}=\sqrt{\ln n / n}$, 
\begin{align}\label{q222}
Q_{22}&\lesssim\sum\limits_{j=j'+1}^{j_{1}}\sum\limits_{u=1}^{2^{d}-1}\sum\limits_{\bk\in\Lambda_{j}}\frac{\ln n}{n}~\bbone_{\{|\beta_{j,\bk,u}|\geq\frac{\kappa t_{n}}{2}\}}\lesssim\sum\limits_{j=j'+1}^{j_{1}}\sum\limits_{u=1}^{2^{d}-1}\sum\limits_{\bk\in\Lambda_{j}}\frac{\ln n}{n}\left( \frac{\beta_{j,\bk,u}}{\kappa t_{n}/2}\right)^{2}\notag\\
&\lesssim\sum\limits_{j=j'+1}^{j_{1}}2^{-2js}\lesssim2^{-2j's}\thicksim n^{-\frac{2s}{2s+d}}.
\end{align}
When $1\leq p<2$ and $s>d/p$, $B_{p,q}^{s}([0,1]^{d})\subseteq B_{2,\infty}^{s+d/2-d/p}([0,1]^{d})$. Then
\begin{align}\label{q223}
Q_{22}&\lesssim\sum\limits_{j=j'+1}^{j_{1}}\sum\limits_{u=1}^{2^{d}-1}\sum\limits_{\bk\in\Lambda_{j}}\frac{\ln n}{n}~\bbone_{\{|\beta_{j,\bk,u}|\geq\frac{\kappa t_{n}}{2}\}}\lesssim\sum\limits_{j=j'+1}^{j_{1}}\sum\limits_{u=1}^{2^{d}-1}\sum\limits_{\bk\in\Lambda_{j}}\frac{\ln n}{n}\left( \frac{\beta_{j,\bk,u}}{\kappa t_{n}/2}\right)^{p}\notag\\
&\lesssim\sum\limits_{j=j'+1}^{j_{1}}(\ln n)n^{\frac{p}{2}-1}2^{-j(s+d/2-d/p)p}\notag\\
&\lesssim (\ln n)n^{\frac{p}{2}-1}2^{-j'(s+d/2-d/p)p}\thicksim (\ln n)n^{-\frac{2s}{2s+d}}.
\end{align}
It follows from the upper bounds above that
\begin{eqnarray}\label{qtt}
Q_{2}\lesssim(\ln n) n^{-\frac{2s}{2s+d}}.
\end{eqnarray}
Finally, one evaluates $Q_{3}$. Clearly,
\begin{align}
Q_{31}&:=\sum\limits_{j=j_{*}}^{j'}\sum\limits_{u=1}^{2^{d}-1}\sum\limits_{\bk\in\Lambda_{j}}\left|\beta_{j,\bk,u}\right|^{2}\bbone_{\{|\beta_{j,\bk,u}|\leq2\kappa t_{n}\}}\notag\\
&\leq\sum\limits_{j=j_{*}}^{j'}\sum\limits_{u=1}^{2^{d}-1}\sum\limits_{\bk\in\Lambda_{j}}\Big|2\kappa t_{n}\Big|^{2}
\lesssim\sum\limits_{j=j_{*}}^{j'}\frac{\ln n}{n}2^{jd}\lesssim\frac{\ln n}{n}2^{j'd}. \notag
\end{align}
This with the choice of $2^{j'}$ shows
\begin{equation*}
Q_{31}\lesssim(\ln n)n^{-\frac{2s}{2s+d}}.
\end{equation*}
On the other hand, $Q_{32}:=\sum\limits_{j=j'+1}^{j_{1}}\sum\limits_{u=1}^{2^{d}-1}\sum\limits_{\bk\in\Lambda_{j}}\left|\beta_{j,\bk,u}\right|^{2}\bbone_{\{|\beta_{j,\bk,u}|\leq2\kappa t_{n}\}}$.
According to the arguments of \eqref{q222}, for $p\geq2$, 
\begin{equation*}
Q_{32}\lesssim\sum\limits_{j=j'+1}^{j_{1}}\sum\limits_{u=1}^{2^{d}-1}\sum\limits_{\bk\in\Lambda_{j}}\left|\beta_{j,\bk,u}\right|^{2}\lesssim n^{-\frac{2s}{2s+d}}.
\end{equation*}
 When $1\leq p<2$, $\left|\beta_{j,\bk,u}\right|^{2}\bbone_{\{|\beta_{j,\bk,u}|\leq2\kappa t_{n}\}}\leq\left|\beta_{j,\bk,u}\right|^{p}\left|2\kappa t_{n}\right|^{2-p}$. Then similar to the arguments of \eqref{q223},
\begin{align*}
Q_{32}&\lesssim\sum\limits_{j=j'+1}^{j_{1}}\sum\limits_{u=1}^{2^{d}-1}\sum\limits_{\bk\in\Lambda_{j}}\left|\beta_{j,\bk,u}\right|^{p}\left|2\kappa t_{n}\right|^{2-p}\notag\\
&\lesssim\left(\frac{\ln n}{n}\right)^{\frac{2-p}{2}}\sum\limits_{j=j'+1}^{j_{1}}2^{-j(s+d/2-d/p)p}\lesssim\left(\frac{\ln n}{n}\right)^{\frac{2-p}{2}}2^{-j'(s+d/2-d/p)p}\notag\\
&\lesssim\left(\frac{\ln n}{n}\right)^{\frac{2-p}{2}}\left(\frac{1}{n}\right)^{\frac{(s+d/2-d/p)p}{2s+d}}\leq (\ln n)n^{-\frac{2s}{2s+d}}.
\end{align*}
It follows from the inequalities above that
\begin{equation}\label{qt}
Q_{3}\lesssim (\ln n)n^{-\frac{2s}{2s+d}},
\end{equation} 
in both cases.
Owing to \eqref{q2}, \eqref{q3}, \eqref{qtt}, and \eqref{qt}, we prove that
\begin{equation*}
Q\lesssim(\ln n)n^{-\frac{2s}{2s+d}},
\end{equation*}
which is the desired conclusion. $\hfill \Box$

\section*{Acknowledgements}
We would like to thank the reviewers for their thoughtful comments
that have improved the manuscript.

\bibliographystyle{chicago} 
\bibliography{reg-multi}

\begin{thebibliography}{}

\bibitem[\protect\citeauthoryear{Abramovich, Bailey, and Sapatinas}{Abramovich
  et~al.}{2000}]{anto2}
Abramovich, F., T.~C. Bailey, and T.~Sapatinas (2000).
\newblock Wavelet analysis and its statistical applications.
\newblock {\em Journal of the Royal Statistical Society: Series D (The
  Statistician)\/}~{\em 49\/}(1), 1--29.

\bibitem[\protect\citeauthoryear{Brown, Levine, et~al.}{Brown
  et~al.}{2007}]{brown}
Brown, L.~D., M.~Levine, et~al. (2007).
\newblock Variance estimation in nonparametric regression via the difference
  sequence method.
\newblock {\em The Annals of Statistics\/}~{\em 35\/}(5), 2219--2232.

\bibitem[\protect\citeauthoryear{Cai and Brown}{Cai and Brown}{1999}]{cai:99}
Cai, T.~T. and L.~D. Brown (1999).
\newblock Wavelet estimation for samples with random uniform design.
\newblock {\em Statistics \& Probability Letters\/}~{\em 42\/}(3), 313--321.

\bibitem[\protect\citeauthoryear{Cai, Brown, et~al.}{Cai et~al.}{1998}]{cai:98}
Cai, T.~T., L.~D. Brown, et~al. (1998).
\newblock Wavelet shrinkage for nonequispaced samples.
\newblock {\em The Annals of Statistics\/}~{\em 26\/}(5), 1783--1799.

\bibitem[\protect\citeauthoryear{Cai, Wang, et~al.}{Cai et~al.}{2008}]{cai}
Cai, T.~T., L.~Wang, et~al. (2008).
\newblock Adaptive variance function estimation in heteroscedastic
  nonparametric regression.
\newblock {\em The Annals of Statistics\/}~{\em 36\/}(5), 2025--2054.

\bibitem[\protect\citeauthoryear{Chaubey, Chesneau, and Doosti}{Chaubey
  et~al.}{2015}]{chaubey}
Chaubey, Y.~P., C.~Chesneau, and H.~Doosti (2015).
\newblock Adaptive wavelet estimation of a density from mixtures under
  multiplicative censoring.
\newblock {\em Statistics\/}~{\em 49\/}(3), 638--659.

\bibitem[\protect\citeauthoryear{Chesneau}{Chesneau}{2013}]{chesneau}
Chesneau, C. (2013).
\newblock On the adaptive wavelet estimation of a multidimensional regression
  function under $\alpha$-mixing dependence: Beyond the standard assumptions on
  the noise.
\newblock {\em Commentationes Mathematicae Universitatis Carolinae\/}~{\em 4},
  527--556.

\bibitem[\protect\citeauthoryear{Chesneau, Kou, and Navarro}{Chesneau
  et~al.}{2019}]{chesneau2018linear}
Chesneau, C., J.~Kou, and F.~Navarro (2019).
\newblock {Linear wavelet estimation in regression with additive and
  multiplicative noise}.
\newblock working paper or preprint.

\bibitem[\protect\citeauthoryear{Chichignoud}{Chichignoud}{2012}]{chi}
Chichignoud, M. (2012).
\newblock Minimax and minimax adaptive estimation in multiplicative regression:
  locally bayesian approach.
\newblock {\em Probability Theory and Related Fields\/}~{\em 153\/}(3-4),
  543--586.

\bibitem[\protect\citeauthoryear{Cohen, Daubechies, and Vial}{Cohen
  et~al.}{1993}]{cohen}
Cohen, A., I.~Daubechies, and P.~Vial (1993).
\newblock Wavelets on the interval and fast wavelet transforms.
\newblock {\em Applied and Computational Harmonic Analysis\/}~{\em 1\/}(1),
  54--81.

\bibitem[\protect\citeauthoryear{Comte}{Comte}{2015}]{comte}
Comte, F. (2015).
\newblock {\em Estimation non-param{\'e}trique}.
\newblock Spartacus-IDH.

\bibitem[\protect\citeauthoryear{Daouia and Simar}{Daouia and
  Simar}{2005}]{daouia2005robust}
Daouia, A. and L.~Simar (2005).
\newblock Robust nonparametric estimators of monotone boundaries.
\newblock {\em Journal of Multivariate Analysis\/}~{\em 96\/}(2), 311--331.

\bibitem[\protect\citeauthoryear{Daubechies}{Daubechies}{1992}]{daub}
Daubechies, I. (1992).
\newblock {\em Ten lectures on wavelets}, Volume~61.
\newblock Siam.

\bibitem[\protect\citeauthoryear{De~Prins, Simar, and Tulkens}{De~Prins
  et~al.}{1984}]{de1984measuring}
De~Prins, D., L.~Simar, and H.~Tulkens (1984).
\newblock Measuring labour efficiency.
\newblock {\em Post Offices in The Performance of Public Enterprises: Concepts
  and Measurement (P. Pestieau ve H. Tulkens M. Marchand)\/}, 243--267.

\bibitem[\protect\citeauthoryear{Delyon and Juditsky}{Delyon and
  Juditsky}{1996}]{delyon}
Delyon, B. and A.~Juditsky (1996).
\newblock On minimax wavelet estimators.
\newblock {\em Applied and Computational Harmonic Analysis\/}~{\em 3\/}(3),
  215--228.

\bibitem[\protect\citeauthoryear{Donoho, Johnstone, Kerkyacharian, and
  Picard}{Donoho et~al.}{1995}]{donoho3}
Donoho, D.~L., I.~M. Johnstone, G.~Kerkyacharian, and D.~Picard (1995).
\newblock Wavelet shrinkage: asymptopia?
\newblock {\em Journal of the Royal Statistical Society. Series B
  (Methodological)\/}, 301--369.

\bibitem[\protect\citeauthoryear{Fan, Li, and Weersink}{Fan
  et~al.}{1996}]{fan1996semiparametric}
Fan, Y., Q.~Li, and A.~Weersink (1996).
\newblock Semiparametric estimation of stochastic production frontier models.
\newblock {\em Journal of Business \& Economic Statistics\/}~{\em 14\/}(4),
  460--468.

\bibitem[\protect\citeauthoryear{Farrell}{Farrell}{1957}]{farrell1957measurement}
Farrell, M.~J. (1957).
\newblock The measurement of productive efficiency.
\newblock {\em Journal of the Royal Statistical Society: Series A
  (General)\/}~{\em 120\/}(3), 253--281.

\bibitem[\protect\citeauthoryear{Gijbels, Mammen, Park, and Simar}{Gijbels
  et~al.}{1999}]{gijbels1999estimation}
Gijbels, I., E.~Mammen, B.~U. Park, and L.~Simar (1999).
\newblock On estimation of monotone and concave frontier functions.
\newblock {\em Journal of the American Statistical Association\/}~{\em
  94\/}(445), 220--228.

\bibitem[\protect\citeauthoryear{Girard, Guillou, and Stupfler}{Girard
  et~al.}{2013}]{girard2013frontier}
Girard, S., A.~Guillou, and G.~Stupfler (2013).
\newblock Frontier estimation with kernel regression on high order moments.
\newblock {\em Journal of Multivariate Analysis\/}~{\em 116}, 172--189.

\bibitem[\protect\citeauthoryear{Girard and Jacob}{Girard and
  Jacob}{2008}]{girard2008frontier}
Girard, S. and P.~Jacob (2008).
\newblock Frontier estimation via kernel regression on high power-transformed
  data.
\newblock {\em Journal of Multivariate Analysis\/}~{\em 99\/}(3), 403--420.

\bibitem[\protect\citeauthoryear{Hall and Marron}{Hall and
  Marron}{1990}]{hall1990variance}
Hall, P. and J.~Marron (1990).
\newblock On variance estimation in nonparametric regression.
\newblock {\em Biometrika\/}~{\em 77\/}(2), 415--419.

\bibitem[\protect\citeauthoryear{Hall, Turlach, et~al.}{Hall
  et~al.}{1997}]{hall:97}
Hall, P., B.~A. Turlach, et~al. (1997).
\newblock Interpolation methods for nonlinear wavelet regression with
  irregularly spaced design.
\newblock {\em The Annals of Statistics\/}~{\em 25\/}(5), 1912--1925.

\bibitem[\protect\citeauthoryear{H{\"a}rdle, Kerkyacharian, Picard, and
  Tsybakov}{H{\"a}rdle et~al.}{2012}]{hardle}
H{\"a}rdle, W., G.~Kerkyacharian, D.~Picard, and A.~Tsybakov (2012).
\newblock {\em Wavelets, approximation, and statistical applications}, Volume
  129.
\newblock Springer Science \& Business Media.

\bibitem[\protect\citeauthoryear{H{\"a}rdle and Tsybakov}{H{\"a}rdle and
  Tsybakov}{1997}]{hardletsy}
H{\"a}rdle, W. and A.~Tsybakov (1997).
\newblock Local polynomial estimators of the volatility function in
  nonparametric autoregression.
\newblock {\em Journal of Econometrics\/}~{\em 81\/}(1), 223--242.

\bibitem[\protect\citeauthoryear{Huang, Pi, and Progri}{Huang
  et~al.}{2013}]{huang}
Huang, P., Y.~Pi, and I.~Progri (2013).
\newblock Gps signal detection under multiplicative and additive noise.
\newblock {\em The Journal of Navigation\/}~{\em 66\/}(4), 479--500.

\bibitem[\protect\citeauthoryear{Jirak, Meister, Rei{\ss}, et~al.}{Jirak
  et~al.}{2014}]{jirak2014adaptive}
Jirak, M., A.~Meister, M.~Rei{\ss}, et~al. (2014).
\newblock Adaptive function estimation in nonparametric regression with
  one-sided errors.
\newblock {\em The Annals of Statistics\/}~{\em 42\/}(5), 1970--2002.

\bibitem[\protect\citeauthoryear{Korostelev and Tsybakov}{Korostelev and
  Tsybakov}{2012}]{korostelev2012minimax}
Korostelev, A.~P. and A.~B. Tsybakov (2012).
\newblock {\em Minimax theory of image reconstruction}, Volume~82.
\newblock Springer Science \& Business Media.

\bibitem[\protect\citeauthoryear{Kuan, Sawchuk, Strand, and Chavel}{Kuan
  et~al.}{1985}]{kuan1985adaptive}
Kuan, D.~T., A.~A. Sawchuk, T.~C. Strand, and P.~Chavel (1985).
\newblock Adaptive noise smoothing filter for images with signal-dependent
  noise.
\newblock {\em IEEE Transactions on Pattern Analysis \& Machine
  Intelligence\/}~(2), 165--177.

\bibitem[\protect\citeauthoryear{Kulik, Raimondo, et~al.}{Kulik
  et~al.}{2009}]{kulik:09}
Kulik, R., M.~Raimondo, et~al. (2009).
\newblock Wavelet regression in random design with heteroscedastic dependent
  errors.
\newblock {\em The Annals of Statistics\/}~{\em 37\/}(6A), 3396--3430.

\bibitem[\protect\citeauthoryear{Kulik, Wichelhaus, et~al.}{Kulik
  et~al.}{2011}]{kulik2011nonparametric}
Kulik, R., C.~Wichelhaus, et~al. (2011).
\newblock Nonparametric conditional variance and error density estimation in
  regression models with dependent errors and predictors.
\newblock {\em Electronic Journal of Statistics\/}~{\em 5}, 856--898.

\bibitem[\protect\citeauthoryear{Kumbhakar, Park, Simar, and Tsionas}{Kumbhakar
  et~al.}{2007}]{kumbhakar2007nonparametric}
Kumbhakar, S.~C., B.~U. Park, L.~Simar, and E.~G. Tsionas (2007).
\newblock Nonparametric stochastic frontiers: a local maximum likelihood
  approach.
\newblock {\em Journal of Econometrics\/}~{\em 137\/}(1), 1--27.

\bibitem[\protect\citeauthoryear{Mallat}{Mallat}{2008}]{mallat:08}
Mallat, S. (2008).
\newblock {\em A wavelet tour of signal processing: the sparse way}.
\newblock Academic press.

\bibitem[\protect\citeauthoryear{Mateo and Fern{\'a}ndez-Caballero}{Mateo and
  Fern{\'a}ndez-Caballero}{2009}]{mateo2009finding}
Mateo, J.~L. and A.~Fern{\'a}ndez-Caballero (2009).
\newblock Finding out general tendencies in speckle noise reduction in
  ultrasound images.
\newblock {\em Expert Systems with Applications\/}~{\em 36\/}(4), 7786--7797.

\bibitem[\protect\citeauthoryear{Meyer}{Meyer}{1992}]{meyer}
Meyer, Y. (1992).
\newblock {\em Wavelets and operators}, Volume~1.
\newblock Cambridge university press.

\bibitem[\protect\citeauthoryear{Muller, Stadtmuller, et~al.}{Muller
  et~al.}{1987}]{muller1987estimation}
Muller, H.-G., U.~Stadtmuller, et~al. (1987).
\newblock Estimation of heteroscedasticity in regression analysis.
\newblock {\em The Annals of Statistics\/}~{\em 15\/}(2), 610--625.

\bibitem[\protect\citeauthoryear{Nason}{Nason}{1996}]{nason:96}
Nason, G.~P. (1996).
\newblock Wavelet shrinkage using cross-validation.
\newblock {\em Journal of the Royal Statistical Society. Series B
  (Methodological)\/}, 463--479.

\bibitem[\protect\citeauthoryear{Navarro and Chesneau}{Navarro and
  Chesneau}{2019}]{chesnav}
Navarro, F. and C.~Chesneau (2019).
\newblock {\em R package {rwavelet}: Wavelet Analysis}.
\newblock (Version 0.4.0).

\bibitem[\protect\citeauthoryear{Navarro and Saumard}{Navarro and
  Saumard}{2017}]{navarro17}
Navarro, F. and A.~Saumard (2017).
\newblock Slope heuristics and v-fold model selection in heteroscedastic
  regression using strongly localized bases.
\newblock {\em ESAIM: Probability and Statistics\/}~{\em 21}, 412--451.

\bibitem[\protect\citeauthoryear{Navarro and Saumard}{Navarro and
  Saumard}{2018}]{navarro2}
Navarro, F. and A.~Saumard (2018).
\newblock Efficiency of the v -fold model selection for localized bases.
\newblock In P.~Bertail, D.~Blanke, P.-A. Cornillon, and E.~Matzner-L{\o}ber
  (Eds.), {\em Nonparametric Statistics}, Cham, pp.\  53--68. Springer
  International Publishing.

\bibitem[\protect\citeauthoryear{Rabbani, Vafadust, Abolmaesumi, and
  Gazor}{Rabbani et~al.}{2008}]{rabbani2008speckle}
Rabbani, H., M.~Vafadust, P.~Abolmaesumi, and S.~Gazor (2008).
\newblock Speckle noise reduction of medical ultrasound images in complex
  wavelet domain using mixture priors.
\newblock {\em IEEE Transactions on Biomedical Engineering\/}~{\em 55\/}(9),
  2152--2160.

\bibitem[\protect\citeauthoryear{Shen, Gao, Witten, and Han}{Shen
  et~al.}{2019}]{shen2019optimal}
Shen, Y., C.~Gao, D.~Witten, and F.~Han (2019).
\newblock Optimal estimation of variance in nonparametric regression with
  random design.
\newblock {\em arXiv preprint arXiv:1902.10822\/}.

\bibitem[\protect\citeauthoryear{Simar and Zelenyuk}{Simar and
  Zelenyuk}{2011}]{simar2011stochastic}
Simar, L. and V.~Zelenyuk (2011).
\newblock Stochastic fdh/dea estimators for frontier analysis.
\newblock {\em Journal of Productivity Analysis\/}~{\em 36\/}(1), 1--20.

\bibitem[\protect\citeauthoryear{Triebel}{Triebel}{1994}]{triebel}
Triebel, H. (1994).
\newblock Theory of function spaces ii.
\newblock {\em Bull. Amer. Math. Soc\/}~{\em 31}, 119--125.

\bibitem[\protect\citeauthoryear{Tsybakov}{Tsybakov}{2009}]{tsybakov}
Tsybakov, A.~B. (2009).
\newblock Introduction to nonparametric estimation. revised and extended from
  the 2004 french original. translated by vladimir zaiats.

\bibitem[\protect\citeauthoryear{Verzelen, Gassiat, et~al.}{Verzelen
  et~al.}{2018}]{verzelen2018adaptive}
Verzelen, N., E.~Gassiat, et~al. (2018).
\newblock Adaptive estimation of high-dimensional signal-to-noise ratios.
\newblock {\em Bernoulli\/}~{\em 24\/}(4B), 3683--3710.

\bibitem[\protect\citeauthoryear{Wang, Brown, Cai, Levine, et~al.}{Wang
  et~al.}{2008}]{wang}
Wang, L., L.~D. Brown, T.~T. Cai, M.~Levine, et~al. (2008).
\newblock Effect of mean on variance function estimation in nonparametric
  regression.
\newblock {\em The Annals of Statistics\/}~{\em 36\/}(2), 646--664.

\end{thebibliography}
\end{document}